\documentclass[twoside, onecolumn, 10pt]{article}
\usepackage{latexsym}
\usepackage{amssymb}
\usepackage{amsmath}
\usepackage{mathrsfs}
\newtheorem{theorem}{Theorem}[section]
\newtheorem{lemma}[theorem]{Lemma}
\newtheorem{corollary}[theorem]{Corollary}

\newtheorem{problem}[theorem]{Problem}

\begin{document}

\thispagestyle{empty}

\def\fl{\noindent}
\def\nb{\nonumber}

\renewcommand{\baselinestretch}{1.0}

\begin{center}
{\LARGE Formulas for calculating the extremal
ranks and inertias of a matrix-valued function subject to matrix equation restrictions}
\end{center}

\begin{center}
{\large Yongge Tian}
\end{center}

\begin{center}
{\footnotesize {\em CEMA, Central University of Finance and Economics, Beijing 100081, China}}
\end{center}

\renewcommand{\thefootnote}{\fnsymbol{footnote}}
\footnotetext{{\it E-mail Address:} yongge.tian@gmail.com}

\noindent {\bf Abstract.} \ {\small
Matrix rank and inertia optimization problems are a class of discontinuous optimization problems
in which the decision variables are matrices running over certain matrix sets, while
the ranks and inertias of the variable matrices are taken as integer-valued objective functions.
 In this paper, we establish a group of explicit formulas for calculating the maximal and minimal
 values of the rank and inertia objective functions of the Hermitian matrix  expression $A_1 - B_1XB_1^{*}$
 subject to the common Hermitian solution of a pair of consistent matrix equations $B_2XB^{*}_2 = A_2$
 and $B_3XB_3^{*} = A_3$, and Hermitian solution of the consistent matrix equation $B_4X= A_4$, respectively.
Many consequences are obtained, in particular, necessary and sufficient conditions
are established for the triple matrix equations $B_1XB^{*}_1 =A_1$, $B_2XB^{*}_2 = A_2$ and $B_3XB^{*}_3 = A_3$
to have a common Hermitian solution, as necessary and sufficient conditions for the two matrix equations
$B_1XB^{*}_1 =A_1$ and $B_4X = A_4$ to have a common Hermitian solution.\\

\medskip

\\

}

\renewcommand{\thefootnote}{\fnsymbol{footnote}}
\footnotetext{}

\section{{\bf Introduction}}

\renewcommand{\theequation}{\thesection.\arabic{equation}}
\setcounter{section}{1}

Throughout this paper,
\begin{enumerate}
\item[]
${\mathbb C}^{m\times n}$ stands for the set of all $m\times n$ complex
matrices; \\[-5mm]

${\mathbb C}_{{\rm H}}^{m}$ stands for the set of all $m\times m$ complex Hermitian matrices; \\[-5mm]

\item[] $A^{T}$,  $A^{*}$,  $r(A)$,  ${\mathscr
R}(A)$ and ${\mathscr N}(A)$ stand for the transpose, conjugate
transpose, rank, range (column space) and null space of a matrix
$A\in {\mathbb C}^{m\times n}$, respectively;\\[-5mm]

\item[] $I_m$ denotes the identity matrix of order $m$;\\[-5mm]

\item[] $[\, A, \, B\,]$ denotes a row block
matrix consisting of $A$ and $B$; \\[-5mm]

\item[] $A > 0$ ($A \geqslant 0$) means that $A$ is Hermitian positive definite (Hermitian positive semi-definite);\\[-5mm]

\item[] two $A, \, B \in {\mathbb C}_{{\rm H}}^{m}$ are said to satisfy
the inequality $A >B$ ($A \geqslant B)$ in the L\"owner partial
ordering if $A - B$ is positive definite (positive semi-definite);\\[-5mm]

\item[] the Moore--Penrose inverse of $A\in {\mathbb C}^{m\times n}$, denoted by
$A^{\dag}$, is defined to be the unique solution $X$ satisfying the
four matrix equations $AXA = A$,  $XAX = X$, $(AX)^{*} = AX$ and
$(XA)^{*} =XA$, which satisfies $AA^{\dag} =A^{\dag}A$ if $A = A^{*}$;\\[-5mm]

\item[] a matrix $X$ is called a Hermitian $g$-inverse of
$A \in {\mathbb C}^{m}_{{\rm H}}$, denoted by $A^{-}$, if it
satisfies both $AXA = A$ and $X = X^{*}$;\\[-5mm]

\item[] $E_A$ and $F_A$ stand for $E_A = I_m - AA^{\dag}$ and $F_A = I_n - A^{\dag}A$.
The ranks of $E_A$ and $F_A$ are given by $r(E_A) = m - r(A)$ and $r(F_A)
= n - r(A)$; \\[-5mm]

\item[] $i_{+}(A)$ and $i_{-}(A)$, usually called the partial inertia of $A \in {\mathbb C}^{m}_{{\rm H}}$,
 are defined to be the numbers of the positive and negative eigenvalues of $A$ counted with multiplicities,
respectively, which satisfy  $r(A) = i_{+}(A) + i_{-}(A)$.
\end{enumerate}

The matrix approximation problem is to approximate optimally, with respect to some criteria,
a matrix by one of the same dimension from a given feasible matrix set. Assume that $A$ is a matrix to
be approximated. Then a conventional statement of general matrix optimization  problems of $A$
from this point of view can be written as
\begin{align}
{\rm minimize}\ \rho(\, A - X \,) \ \  {\rm subject \ to}  \ X  \in {\cal S},
\label{ww11}
\end{align}
where $\rho(\cdot)$ is certain objective function, which is usually taken as the determinant, trace,
 norms, rank, inertia of  matrix, and ${\cal S}$  is a given feasible matrix set.
 A best-known case of (\ref{ww11})  is to  minimize the norm $ \|\, A - X \,\|^2_{F}$ subject to  $X  \in {\cal S}$.

In this paper, we take the matrix set  ${\cal S}$ as
\begin{align}
& {\cal S}=  \{ \,\phi(X) =  A_1-B_1XB_1^{*}   \ | \  [\,B_2XB^{*}_2, \, B_3XB^{*}_3 \,] = [\,A_2, \, A_3\,] \, \},
\label{yy11}
\\
& {\cal S}=  \{\, \phi(X) =  A_1-B_1XB_1^{*}  \ | \  B_4X = A_4 \, \},
\label{yy12}
\end{align}
where $A_i \in \mathbb C_{{\rm H}}^{m_i}$ and $B_i \in \mathbb
C^{m_i \times n}$,  $A_4, \ B_4 \in {\mathbb C}^{m \times n}$  are given,
$ i =1,\, 2, \,3$, and $X \in {\mathbb C}_{{\rm H}}^{n}$ is a
variable matrix, and study the following constrained optimization problems:

\begin{problem} \label{TS11}
{\rm  For the constrained linear matrix-valued function in {\rm (\ref{yy11})}$,$
establish explicit formulas for calculating
\begin{align}
& \max_{X \in {\mathbb C}_{{\rm H}}^{n}}\! r(\, A_1 - B_1XB^{*}_1 \,) \ \ \ \ \ \ {\rm s.t.} \ \
[\,B_2XB^{*}_2, \, B_3XB^{*}_3 \,] = [\,A_2, \, A_3\,],
\label{yy13}
\\
&\min_{X \in {\mathbb C}_{{\rm H}}^{n}}\! r(\, A_1 - B_1XB^{*}_1 \,) \ \ \ \ \ \ {\rm s.t.}
\ \ [\,B_2XB^{*}_2, \, B_3XB^{*}_3 \,] = [\,A_2, \, A_3\,],
\label{yy14}
\\
&  \min_{X \in {\mathbb C}_{{\rm H}}^{n}}\! i_{\pm}(\, A_1 - B_1XB^{*}_1 \,) \ \ \ \ \ \ {\rm s.t.} \ \ [\,B_2XB^{*}_2, \, B_3XB^{*}_3 \,] = [\,A_2, \, A_3\,],
\label{yy15}
\\
& \min_{X \in {\mathbb C}_{{\rm H}}^{n}}\! i_{\pm}(\, A_1 - B_1XB^{*}_1 \,) \ \ \ \ \ \ {\rm s.t.} \ \ [\,B_2XB^{*}_2, \, B_3XB^{*}_3 \,] = [\,A_2, \, A_3\,].
\label{yy16}
\end{align}
}
\end{problem}

\begin{problem} \label{TS11a}
{\rm
For the constrained linear matrix-valued function in {\rm (\ref{yy11})}$,$
\begin{enumerate}
\item[{\rm (i)}]  establish necessary and sufficient conditions for the following
three matrix equations
\begin{align}
& [\, B_1XB^{*}_1, \,  B_2XB^{*}_2, \, B_3XB^{*}_3 \,] = [\,A_1, \, A_2, \, A_3\,]
\label{yy17}
\end{align}
to have a common Hermitian solution;

\item[{\rm (ii)}]  establish necessary and sufficient conditions for
$A_1 - B_1XB^{*}_1 > \,(\geqslant, \, <, \, \leqslant \,) \, 0$ to hold for  an
$X \in {\mathbb C}_{{\rm H}}^{n}$  satisfying $[\,B_2XB^{*}_2, \, B_3XB^{*}_3 \,] = [\,A_2, \, A_3\,]$;

\item[{\rm (iii)}]  establish  necessary and sufficient conditions
for $A_1 - B_1XB^{*}_1 > \,(\geqslant, \, <, \, \leqslant \,) \, 0$ to hold for all
$X \in {\mathbb C}_{{\rm H}}^{n}$ satisfying $[\,B_2XB^{*}_2, \, B_3XB^{*}_3 \,] = [\,A_2, \, A_3\,],$
namely$,$ $A_1 - B_1XB^{*}_1$ is a positive map under the restriction $[\,B_2XB^{*}_2, \, B_3XB^{*}_3 \,] = [\,A_2, \, A_3\,].$
\end{enumerate}
}
\end{problem}

\begin{problem} \label{TS11b}
{\rm
For the constrained linear matrix-valued function in {\rm (\ref{yy12})}, establish explicit formulas for calculating
\begin{align}
& \max_{X \in {\mathbb C}_{{\rm H}}^{n}}\! r(\, A_1 - B_1XB^{*}_1 \,) \ \ \ \ \ \ {\rm s.t.} \ \ B_4X = A_4,
\label{yy18}
\\
& \min_{X \in {\mathbb C}_{{\rm H}}^{n}}\! r(\, A_1 - B_1XB^{*}_1 \,)  \ \ \ \ \ \ {\rm s.t.} \ \ B_4X = A_4,
\label{yy19}
\\
& \max_{X \in {\mathbb C}_{{\rm H}}^{n}}\! i_{\pm}(\, A_1 - B_1XB^{*}_1 \,) \ \ \ \ {\rm s.t.} \ \ B_4X  =A_4,
\label{yy110}
\\
& \min_{X \in {\mathbb C}_{{\rm H}}^{n}} \!i_{\pm}(\, A_1 - B_1XB^{*}_1 \,) \ \ \ \ {\rm s.t.} \ \ B_4X = A_4.
\label{yy111}
\end{align}
}
\end{problem}

\begin{problem} \label{TS11c}
{\rm
For the linear matrix  map  in {\rm (\ref{yy12})},
\begin{enumerate}
\item[{\rm (i)}]  establish necessary and sufficient conditions for the following
two matrix equations $[\, B_1XB^{*}_1, \,  B_4X \,] = [\,A_1, \, A_4 \,]$ to have a common
 Hermitian solution and nonnegative definite solution;

\item[{\rm (ii)}]  establish  necessary and sufficient conditions for
$ A_1 - B_1XB^{*}_1 > \,(\geqslant, \, <, \, \leqslant \,) \, 0$ to hold for an
$X \in {\mathbb C}_{{\rm H}}^{n}$ satisfying $B_4X = A_4$;

\item[{\rm (iii)}]  establish  necessary and sufficient conditions for
$A_1 - B_1XB^{*}_1 > \,(\geqslant, \, <, \, \leqslant \,) \, 0$ to hold  for all
$X \in {\mathbb C}_{{\rm H}}^{n}$  satisfying  $B_4X = A_4,$ namely$,$
$A_1 - B_1XB^{*}_1$ is a positive map under the restriction $B_4X = A_4$.
\end{enumerate}
}
\end{problem}

The extremal ranks and inertias of a matrix expression can directly be used to
describe some  behaviors of the matrix expression, for example,

(I) the maximal and minimal dimensions of the row
and column spaces of the matrix expression;

(II)  nonsingularity of the matrix expression when it is square;

(III) solvability of the corresponding matrix equation;

(IV) rank, inertia and range invariance of the matrix expression;

(V) semi-definiteness of a matrix expression; etc.\\
On the other hand, matrix rank and inertia optimization problems are NP-hard in general due to the
discontinuity and combinational nature of rank and inertia of a matrix and the complexity of
algebraic structure of ${\cal S}$.

Mappings between matrix spaces with symmetric patterns can be constructed arbitrarily,
but the linear map $\phi(X)$ in (\ref{yy11}) and (\ref{yy12}) is the simplest cases among all
LMFs with symmetric patterns. This matrix-valued function is the starting point in dealing with
 various complicated matrix-valued functions with symmetric patterns.
In recent years, the matrix-valued function $\phi(X) = A - BXB^{*}$ was reconsidered and many new results
on its algebraic properties were obtained, for instance,
\begin{enumerate}
\item[{\rm (i)}] Expansion formulas for calculating the (global extremal) rank and inertia of
$\phi(X)$ when $X$ running over ${\mathbb C}_{{\rm H}}^{n}$,  see \cite{LT-nla,T-laa10,TL}.  \\[-4mm]

\item[{\rm (ii)}] Nonsingularity, positive definiteness, rank and inertia invariance, etc., of the
$\phi(X)$,  see \cite{T-laa10,TL}. \\[-4mm]

\item[{\rm (iii)}]  Canonical forms of the $\phi(X)$ under generalized singular value decompositions
and their algebraic properties, see  \cite{LT-nla}. \\[-4mm]

\item[{\rm (iv)}]  Solutions and least-squares solutions of the matrix equation $\phi(X) = 0$ and their algebraic properties,
 see \cite{LT2,LTT,T-laa11,T-med,T-nla}.
\\[-4mm]

\item[{\rm (v)}]  Solutions of the  matrix inequalities  $\phi(X) > \, (\geqslant, \, <, \, \leqslant)\, 0$ and
their properties, see  \cite{T-laa10}. \\[-4mm]

\item[{\rm (vi)}] Minimization of ${\rm tr}[\,\phi(X) \phi^{*}(X)\,]$
s.t. $r[\phi(X)]= \min$, see \cite{T-nla}. \\[-4mm]

\item[{\rm (vii)}]  Formulas for calculating the extremal rank and inertia of
$\phi(X)$ under the restrictions $r(X) \leqslant k$ and/or $\pm X \geqslant 0$, see \cite{T-mcm}. \\[-4mm]

\item[{\rm (viii)}]  Formulas for calculating the extremal rank and inertia of
$\phi(X)$ subject to a consistent matrix equation $CXC^{*} = D$, see \cite{LT4}.
\end{enumerate}
This basic work was also extended to some general LMFs, such as,
 $A - BX - (BX)^{*}$, $A - BXB^{*} - CYC^{*}$ and $A - BXC - (BXC)^{*}$, where $X$ and $Y$ are
 (Hermitian) variable matrices of appropriate sizes; see \cite{CHW,LT1,LT2,LT3,LT4,T-ela11,T-laa11,T-med}.

We shall use some pure algebraic operations on matrices to derive
two groups of  analytical formulas for calculating the global extremal values of
the objective functions in (\ref{yy13})--(\ref{yy16}) and  (\ref{yy18})--(\ref{yy111}), and then to
present a variety of valuable consequences of these formulas.

Since variable entries in a matrix-valued function are often
regarded as continuous variables in some constrained sets, while the
objective functions---the rank and inertia of the matrix-valued function take
values only from a finite set of nonnegative integers,  Hence,
(\ref{yy13})--(\ref{yy16}) and  (\ref{yy18})--(\ref{yy111}) can be regarded as continuous-integer
optimization problems subject to equality constraints. This kind of
non-smooth optimization problems cannot be solved by using various
optimization methods for solving continuous or discrete cases.
There is no rigorous mathematical theory for solving a general rank
and inertia optimization problem due to the discontinuity and nonconvexity of
rank and inertia of matrix. In fact, it has been realized that rank and inertia optimization problems have
deep connections with computational complexity, and are regarded as NP-hard
in general; see, e.g., \cite{CR,FHB1,FHB2,Gee,HKY,HT1,HT2,Lau,Mah,Na,RFP}.
Fortunately, some special rank and inertia optimization problems now can be solved by
pure algebraical methods. In particular, analytical solutions to the rank and inertia
optimization problems of the  $\phi(X)$ in (\ref{yy11}) and (\ref{yy12}), as well as
(\ref{yy13})--(\ref{yy16}) and  (\ref{yy18})--(\ref{yy111})
can be derived algebraically by using generalized inverses of matrices.

The following are some known results for ranks and inertias of matrices and their usefulness, which will be
used in the latter part of this paper.

\begin{lemma} [\cite{T-laa10}] \label{T12}
Let ${\cal H}$ be a matrix set in ${\mathbb C}_{{\rm H}}^{m}.$ Then$,$
\begin{enumerate}
\item[{\rm (a)}] ${\cal H}$ has a matrix $X > 0$  $(X < 0)$ if and only if
$\max_{X\in {\cal H}} i_{+}(X) = m  \ \left(\max_{X\in
{\cal H}} i_{-}(X) = m \right).$

\item[{\rm (b)}] All $X\in {\cal H}$ satisfy $X >0$ $(X < 0),$
namely$,$  ${\cal H}$ is a subset of the cone of  positive definite matrices
(negative definite matrices)$,$ if and only if
$\min_{X\in {\cal H}} i_{+}(X) = m \ \left( \min_{X\in
{\cal H}} i_{-}(X) = m \right).$

\item[{\rm (c)}] ${\cal H}$ has a matrix  $X \geqslant 0$ $(X \leqslant 0)$ if and only if
$\min_{X \in {\cal H}} i_{-}(X) = 0 \ \left( \min_{X\in
{\cal H}} i_{+}(X) = 0 \right).$

\item[{\rm (d)}] All $X\in {\cal H}$ satisfy $X \geqslant 0$ $(X \leqslant 0)$
namely$,$  ${\cal H}$ is a subset of the cone of nonnegative definite matrices
(semi-definite matrices)$,$ if and only if $\max_{X\in {\cal H}} i_{-}(X) = 0
\ \left(\max_{X\in {\cal H}} i_{+}(\,X) = 0 \right).$
\end{enumerate}
\end{lemma}

The question of whether a given function is  (definite or semi-definite everywhere is ubiquitous in mathematics and
applications. Lemma \ref{T12}(a)--(d) show that if some explicit formulas for calculating the
global maximal and minimal inertias of a given Hermitian matrix-valued function are established,
we can use them, as demonstrated in Sections below, to derive necessary
and sufficient conditions for the Hermitian matrix-valued function to be definite or semi-definite.

\begin{lemma}[\cite{MS}]\label{T13}
Let $A \in \mathbb C^{m \times n}, \, B \in \mathbb C^{m \times p}$ and
$C \in \mathbb C^{q \times n}.$ Then$,$ the following rank expansion formulas hold
\begin{align}
r[\, A, \, B \,] & = r(A) + r(E_AB) = r(B) + r(E_BA),
\label{1.8}
\\
r \!\left[\!\!\begin{array}{c} A \\ C
\end{array}
\!\!\right] & = r(A) + r(CF_A) = r(C) + r(AF_C),
 \label{1.9}
\\
r\!\left[\!\! \begin{array}{cc}  A  & B  \\ C & 0 \end{array}
\!\!\right] & =  r(B) + r(C) +
 r(E_BAF_C).
 \label{1.10}
\end{align}
\end{lemma}

Three useful rank expansion formulas derived from (\ref{1.10}) are
\begin{align}
r\!\left[\!\! \begin{array}{ccc} A & B & 0  \\ C & 0 & P
\end{array} \!\!\right] & = r(P) + r\!\left[\!\!\begin{array}{cc}
A & B\\ E_PC & 0 \end{array} \!\!\right]\!,
\label{17}
\\
r\!\left[\! \begin{array}{ccc} A & B \\ C & 0 \\ 0 & Q
\end{array} \!\right] & = r(Q) + r\!\left[\!\!\begin{array}{cc}
A & BF_Q \\ C & 0 \end{array} \!\!\right]\!,
\label{18}
\\
r\!\left[\!\! \begin{array}{ccc} A & B & 0  \\ C & 0 & P \\ 0 & Q & 0
\end{array} \!\!\right] & = r(P) + r(Q) + r\!\left[\!\!\begin{array}{cc}
A & BF_Q \\ E_PC & 0 \end{array} \!\!\right]\!.
\label{19}
\end{align}
We shall use them in Section 2 to simplify ranks of block matrices involving
 $E_P$ and $F_Q$.

\begin{lemma}[\cite{T-laa10}] \label{T13a}
Let $A \in {\mathbb C}_{{\rm H}}^{m},$ $B \in \mathbb
C^{m\times n},$  $D \in {\mathbb C}_{{\rm H}}^{n},$  and let
$$
U = \left[\!\!\begin{array}{cc}  A  & B  \\ B^{*}  & 0 \end{array}
\!\!\right]\!, \ \ V = \left[\!\!\begin{array}{cc}  A  & B  \\
B^{*}  & D \end{array} \!\!\right]\!.
$$
Then$,$  the following  expansion formulas hold
\begin{align}
i_{\pm}(U)& = r(B) + i_{\pm}(E_BAE_B),
 \label{18j}
\\
i_{\pm}(V) & =i_{\pm}(A) +
 i_{\pm}\!\left[\!\!\begin{array}{cc} 0 &  E_AB
 \\
 B^{*}E_A & D - B^{*}A^{\dag}B \end{array}\!\!\right]\!.
\label{19j}
\end{align}
\begin{enumerate}
\item[{\rm(a)}] If $A \geqslant 0,$ then
\begin{align}
i_{+}(U)= r[\, A,  \, B \,], \ \ i_{-}(U) = r(B), \ \
 r(U) = r[\, A,  \, B \,] + r(B).
\label{111j}
\end{align}

\item[{\rm(b)}] If $A \leqslant 0,$ then
\begin{align}
i_{+}(U) = r(B), \ \ i_{-}(U) = r[\, A,  \, B \,], \ \ r(U) = r[\,
A, \, B \,] + r(B).
\label{112j}
\end{align}

\item[{\rm(c)}] If ${\mathscr R}(B) \subseteq {\mathscr R}(A),$ then
\begin{align}
i_{\pm}(V) = i_{\pm}(A) + i_{\pm}(\,  D - B^{*}A^{\dag}B \,), \ \ r(V)
= r(A) + r(\,  D - B^{*}A^{\dag}B \,).
\label{113j}
\end{align}

\item[{\rm(d)}] If ${\mathscr R}(B) \cap {\mathscr R}(A) =\{0\}$ and
${\mathscr R}(B^{*}) \cap {\mathscr R}(D) =\{0\},$ then
\begin{align}
i_{\pm}(V) = i_{\pm}(A) + i_{\pm}(D) + r(B), \ \ r(V) = r(A) + 2r(B)
+ r(D).
\label{114j}
\end{align}
\end{enumerate}
\end{lemma}

Three general expansion formulas derived from (\ref{18j}) are
\begin{align}
i_{\pm}\!\left[\!\! \begin{array}{cc} A
 & BF_P \\  F_P B^{*}  & 0 \end{array} \!\!\right]  = i_{\pm}\!\left[\!\! \begin{array}{ccc} A
 & B & 0 \\  B^{*} & 0 & P^{*} \\  0 & P & 0 \end{array} \!\!\right] - r(P),
\  r\!\left[\!\! \begin{array}{cc} A
 & BF_P \\  F_P B^{*}  & 0 \end{array} \!\!\right]  = r\!\left[\!\! \begin{array}{ccc} A
 & B & 0 \\  B^{*} & 0 & P^{*} \\  0 & P & 0 \end{array} \!\!\right] - 2r(P).
\label{114u}
\end{align}
We shall use them to simplify the inertias of block Hermitian
matrices that involve $F_P = I - P^{\dag}P$.

\begin{lemma} \label{T13b}
Let $A_j \in \mathbb C^{m_j\times n},$ $B_j\in \mathbb C^{p\times q_j}$ and $C_j\in \mathbb C^{m_j\times q_j}$
 be given$,$ $j =1,2.$ Then$,$
\begin{enumerate}
\item[{\rm(a)}] {\rm \cite{OA}} The pair of matrix equations
\begin{equation}
A_1XB_1 = C_1 \ \  and  \ \ A_2XB_2 = C_2
\label{113}
\end{equation}
have a common solution for $X\in \mathbb C^{n\times p}$ if and only if
\begin{equation}
{\mathscr R}(C_j) \subseteq {\mathscr R}(A_j),  \ \ {\mathscr R}(C^{*}_j) \subseteq {\mathscr R}(B^{*}_j), \ \
r\!\left[\!\begin{array}{ccc} C_1 & 0 & A_1  \\ 0 & -C_2  & A_2 \\ B_1 & B_2 & 0 \end{array}\!\!\right]
=r\!\left[\!\begin{array}{c} A_1  \\ A_2 \end{array}\!\!\right] + r[\, B_1, \, B_2\,], \ \ \ j =1, \, 2.
\label{114}
\end{equation}

\item[{\rm(b)}] {\rm \cite{T-lama00}}  Under {\rm (\ref{118})}$,$ the general common solution
 to {\rm (\ref{117})} can be written in the following parametric form
\begin{equation}
X = X_0 + F_AV_1 +  V_2E_B +  F_{A_1}V_3E_{B_2} +  F_{A_2}V_4E_{B_1},
\label{115}
\end{equation}
where $A =\left[\!\begin{array}{c} A_1  \\ A_2 \end{array}\!\!\right],$
$B =[\, B_1, \, B_2\,],$  and the four matrices $V_1, \ldots, V_4 \in \mathbb C^{n \times p}$
are arbitrary$.$
\end{enumerate}
\end{lemma}

\begin{lemma} \label{T14}
Let $A \in \mathbb C^{m \times n}$ and $B \in \mathbb C_{{\rm H}}^{m}$ be given$.$
Then$,$
\begin{enumerate}
\item[{\rm(a)}] {\rm \cite{Gr1,KM}} The matrix equation $AXA^{*} = B$ has a solution
$X \in \mathbb C_{{\rm H}}^{n}$ if and only if ${\mathscr R}(B) \subseteq {\mathscr R}(A),$
or equivalently$,$ $AA^{\dag}B = B.$

\item[{\rm(b)}] {\rm \cite{T-laa10}} Under $AA^{\dag}B = B,$ the general Hermitian solution
of $AXA^{*} = B$ can be written in the following two forms
\begin{align}
X & = A^{\dag}B(A^{\dag})^{*} +  U - A^{\dag}AUA^{\dag}A,
\label{15d}
\\
X & = A^{\dag}B(A^{\dag})^{*} +  F_AV + V^{*}F_A,
\label{16d}
\end{align}
where  $U \in \mathbb C_{{\rm H}}^{n}$ and $ V \in \mathbb C^{n \times n}$ are arbitrary$.$

\end{enumerate}
\end{lemma}

More results on properties of solutions of $AXA^{*} = B$ can be found in \cite{LT2,LTT}.

\begin{lemma} [\cite{KM}]\label{T15b}
Let $A, \, B\in \mathbb C^{m\times n}$ be given$.$ Then$,$
\begin{enumerate}
\item[{\rm(a)}]  The matrix equation $AX = B$ has a Hermitian solution
$X \in \mathbb C_{{\rm H}}^{n}$  if and only if ${\mathscr R}(B) \subseteq
{\mathscr R}(A)$ and $AB^{*} =BA^{*}.$ In this case$,$ the general
Hermitian solution of $AX = B$ can be written as
\begin{equation}
X = A^{\dag}B + (A^{\dag}B)^{*} - A^{\dag}BA^{\dag}A + F_AUF_A,
\label{qq121}
\end{equation}
where $U\in \mathbb C_{{\rm H}}^{n}$ is arbitrary$.$

\item[{\rm(b)}] The matrix equation $AX = B$
has a solution $0 \leqslant X \in {\mathbb C}_{{\rm H}}^{n}$ if and only if ${\mathscr R}(B) \subseteq {\mathscr
R}(A),$ $AB^{*} \geqslant 0$ and $r(AB^{*}) =r(B).$ In this case$,$ the
general solution of $AX = B$ can be written as
\begin{equation}
X = B^{*}(AB^{*})^{\dag}B   + F_AUF_A,
\label{qq122}
\end{equation}
where $0 \leqslant U\in {\mathbb C}_{{\rm H}}^{n}$ is arbitrary$.$
\end{enumerate}
\end{lemma}

\begin{lemma}  \label{T15a}
 Let $A \in {\mathbb C}^m_{{\rm H}}$ and $B\in \mathbb C^{m\times n}$ be given$.$ Then$,$
\begin{enumerate}
\item[{\rm(a)}]  {\rm \cite{T-laa10,TL}} The global maximal and minimal ranks and inertias of $A- BXB^{*}$ subject to
$X \in {\mathbb C}^n_{{\rm H}}$  are given by
  \begin{align}
\max_{X\in {\mathbb C}^n_{{\rm H}}} r(\,A-BXB^{*}\,) & = r[\, A, \,B \,],
 \label{qq133}
\\
\min_{X \in {\mathbb C}^n_{{\rm H}}} r(\,A-BXB^{*}\,) & = 2r[\, A, \, B\,]- r\!\left[\!\! \begin{array}{cc}  A & B \\
B^{*} & 0\end{array} \!\!\right]\!,
\label{qq134}
 \\
\max_{X\in {\mathbb C}^n_{{\rm H}}}i_{\pm}(\,A-BXB^{*}\,)& = i_{\pm}\!\left[\!\! \begin{array}{cc}  A & B \\
B^{*} & 0\end{array} \!\!\right]\!,
\label{qq135}
 \\
\min_{X \in {\mathbb C}^n_{{\rm H}}}  i_{\pm}(\,A-BXB^{*}\,) & =  r[\, A, \, B \,]  -
 i_{\mp}\!\left[\!\! \begin{array}{cc}  A & B \\
B^{*} & 0\end{array} \!\!\right]\!.
 \label{qq136}
\end{align}

\item[{\rm(b)}]  {\rm \cite{T-mcm}} The global maximal and minimal ranks and inertias of $A- BXB^{*}$ subject to
$0 \leqslant X \in {\mathbb C}^n_{{\rm H}}$  are given by
\begin{align}
&\max_{0 \leqslant X\in {\mathbb C}^{n}_{{\rm H}}}\!\!\!r(\, A + BXB^{*} \,)
= r[\, A, \, B\,], \  \min_{0 \leqslant X\in {\mathbb C}^{n}_{{\rm H}}}\!\!\!r(\, A + BXB^{*} \,)
 = i_{+}(A) + r[\, A, \, B\,] - i_{+}(M),
\label{qq137}
\\
&\max_{0 \leqslant X\in {\mathbb C}^{n}_{{\rm H}}}\!\!\!i_{+}(\, A + BXB^{*} \,)  = i_{+}(M), \
  \min_{0 \leqslant X\in {\mathbb C}^{n}_{{\rm H}}}\!\!\!i_{+}(\, A + BXB^{*} \,)  = i_{+}(A),
\label{qq138}
\\
&\max_{0 \leqslant X\in {\mathbb C}^{n}_{{\rm H}}}\!\!\!i_{-}(\, A + BXB^{*} \,)  = i_{-}(A), \ \
\min_{0 \leqslant X\in {\mathbb C}^{n}_{{\rm H}}}\!\!\!i_{-}(\, A + BXB^{*} \,)  = r[\, A, \, B\,] - i_{+}(M),
 \label{qq139}
\\
&\max_{0 \leqslant X\in {\mathbb C}^{n}_{{\rm H}}}\!\!\!r(\, A - BXB^{*} \,)
= r[\, A, \, B\,],  \  \min_{0 \leqslant X\in {\mathbb C}^{n}_{{\rm H}}}\!\!\!r(\, A - BXB^{*} \,) =
i_{-}(A) + r[\, A, \, B\,] - i_{-}(M),
\label{qq140}
\\
&\max_{0 \leqslant X\in {\mathbb C}^{n}_{{\rm H}}}\!\!\!i_{+}(\, A - BXB^{*} \,)  = i_{+}(A), \ \
\min_{0 \leqslant X\in {\mathbb C}^{n}_{{\rm H}}}\!\!\!i_{+}(\, A - BXB^{*} \,) =  r[\,A, \, B\,] -  i_{-}(M),
\label{qq141}
\\
&\max_{0 \leqslant X\in {\mathbb C}^{n}_{{\rm H}}}\!\!\!i_{-}(\, A - BXB^{*} \,)  = i_{-}(M), \ \
\min_{0 \leqslant X\in {\mathbb C}^{n}_{{\rm H}}}\!\!\!i_{-}(\, A - BXB^{*} \,) = i_{-}(A).
 \label{qq142}
\end{align}
\end{enumerate}
\end{lemma}

\begin{lemma}  [\cite{LT4}] \label{T18}
Let $A \in \mathbb C_{{\rm H}}^{m}$, $B \in \mathbb C^{m\times p}$
and $C \in \mathbb C^{q \times m}$ be given$,$ and let
\begin{align}
& M_1 =\left[\!\!\begin{array}{ccc} A  & B  \\  B^{*}  & 0
\end{array}\!\!\right]\!, \ \ M_2 =\left[\!\!\begin{array}{ccc} A  &
C^{*}  \\  C  & 0   \end{array}\!\!\right]\!, \ \
\label{116}
\\
& N = [\,A,\, B,\, C^{*}\,], \ \   N_1 =
\left[\!\!\begin{array}{ccc} A & B & C^{*} \\ B^{*} & 0 & 0
\end{array}\!\!\right]\!, \ \ N_2 =\left[\!\!\begin{array}{ccc} A & B  & C^{*} \\ C & 0 & 0
\end{array}\!\!\right]\!.
\label{117}
\end{align}
Then$,$ the global maximal and minimal ranks and partial inertias of $A-BXC-(BXC)^{*}$ are given by
\begin{align}
\max_{X\in \mathbb C^{p \times q}}\!\!r[\,A-BXC-(BXC)^{*}\,] & = \min
\left\{ r(N), \ \ r(N_1), \ \ r(N_2) \right\}\!,
\label{118}
\\
\min_{X\in \mathbb C^{p \times q}}\!\!r[\,A-BXC-(BXC)^{*}\,] & =
2r(N) + \max\{\, s_1, \ \ s_2, \ \ s_3, \ \ s_4 \, \},
\label{119}
\\
\max_{X\in \mathbb C^{p \times q}}\!\!i_{\pm}[\,A-BXC-(BXC)^{*}\,] & =
\min\!\left\{i_{\pm}(M_1), \ \ i_{\pm}(M_2) \right\}\!,
\label{120}
\\
\min_{X\in \mathbb C^{p \times q}}\!\!i_{\pm}[\,A-BXC-(BXC)^{*}\,]& =
r(N) + \max\{\, i_{\pm}(M_1) - r(N_1), \ \ i_{\pm}(M_2) - r(N_2)\,
\},
 \label{121}
\end{align}
where
$$
s_1 = r(M_1) - 2r(N_1),  \ \ s_2 = r(M_2) - 2r(N_2), \ \
$$
$$
s_3 = i_{+}(M_1) + i_{-}(M_2) - r(N_1) -  r(N_2), \ \
s_4 = i_{+}(M_1) + i_{-}(M_2) - r(N_1) -  r(N_2).
$$
In particular$,$ if ${\mathscr R}(C^{*})\subseteq {\mathscr R}(B),$ then
\begin{align}
\max_{X\in \mathbb C^{p \times q}}r[\, A - BXC - (BXC)^{*}\,] & =
\min \left\{ r[\,A,\,  B\,],
 \ \ r\!\left[\!\!\begin{array}{cc} A & C^{*}
\\ C & 0
\end{array}\!\!\right] \right\}\!,
\label{122}
\\
\min_{X\in \mathbb C^{p \times q}}\!\!r[\,A-BXC-(BXC)^{*}\,] & =
2r[\,A,\, B\,]+ r\!\left[\!\!\begin{array}{cc} A & C^{*} \\ C & 0
 \end{array}\!\!\right]- 2r\!\left[\!\!\begin{array}{ccc} A & B \\ C & 0
 \end{array}\!\!\right]\!,
\label{123}
\\
\max_{X\in \mathbb C^{p \times q}}\!\!i_{\pm}[\,A-BXC-(BXC)^{*}\,] & =
i_{\pm}\!\left[\!\!\begin{array}{ccc} A
 & C^{*}  \\  C  & 0
   \end{array}\!\!\right]\!,
\label{124}
\\
\min_{X\in \mathbb C^{p \times q}}\!\!i_{\pm}[\,A-BXC-(BXC)^{*}\,] & =
r[\,A,\, B\,] + i_{\pm}\!\left[\!\!\begin{array}{cc} A & C^{*} \\ C
& 0\end{array}\!\!\right]
 - r\!\left[\!\!\begin{array}{ccc} A & B \\ C & 0
\end{array}\!\!\right]\!.
\label{125}
\end{align}
\end{lemma}

The matrices $X$ that satisfy (\ref{118})--(\ref{121}) (namely, the
global maximizers and minimizers of the objective rank and inertia
functions) are not necessarily unique and their expressions were
also given in \cite{LT4} by using certain simultaneous decomposition
of the three given matrices. Observe that the right-hand sides of
(\ref{118})--(\ref{121}) are represented in analytical forms of the
ranks and inertias of the five given block matrices, we can easily
use them to derive extremal ranks and inertias of some general
linear and  nonlinear matrix-valued functions. In these cases, combining
the rank and inertia formulas obtained with the assertions in Lemma
1.1 may yield various conclusions on algebraic properties of linear
and nonlinear matrix-valued functions.

\section{The extremal ranks and inertias of $A - B_1XB^{*}_1$ subject to $B_2XB^{*}_2 = A_2$ and  $B_3 XB^{*}_3 = A_3$}

\renewcommand{\theequation}{\thesection.\arabic{equation}}
\setcounter{section}{3} \setcounter{equation}{0}

We first derive a parametric form for the general common Hermitian solution of
the pair of matrix equations in (\ref{yy11}).

\begin{lemma}[\cite{T-laa11}] \label{T31}
Let $A_i \in \mathbb C_{{\rm H}}^{m_i},$ $B_i \in \mathbb C^{m_i
\times n}$ be given for $ i = 2,\, 3,$ and suppose that
 each of the two matrix equations
\begin{align}
 B_2XB^{*}_2 = A_2 \ \ and  \ \ B_3 XB^{*}_3 = A_3
\label{30}
\end{align}
has a solution$,$ i.e.$,$ ${\mathscr R}(A_i) \subseteq {\mathscr R}(B_i)$  for $i =
2,\, 3.$  Then$,$
\begin{enumerate}
\item[{\rm(a)}] The pair of matrix equations have a common Hermitian solution if and only if
\begin{align}
r\!\left[\!\! \begin{array}{ccc}  A_2  & 0  & B_2  \\ 0 &  -A_3  & B_3  \\
B^{*}_2 & B^{*}_3 & 0 \end{array} \!\! \right]  = 2r\!\left[\!\!
\begin{array}{c} B_2  \\ B_3 \end{array} \!\! \right].
\label{31}
\end{align}

\item[{\rm(b)}] Under {\rm (\ref{31})}$,$ the general common Hermitian solution of
the pair of equations can be written  in the following parametric form
\begin{align}
 X = X_0 + V F_B  +  F_B V^{*} + F_{B_2}U F_{B_3} +
F_{B_3}U^{*}F_{B_2},
\label{32}
\end{align}
 where $X_0$ is a special Hermitian common solution to the pair of equations$,$
 $B = \left[\!\!\begin{array}{c}
B_2  \\ B_3\end{array} \!\!\right]\!,$ and $U, \, V \in \mathbb
C^{n \times n}$ are arbitrary$.$
\end{enumerate}
\end{lemma}

Substituting (\ref{32}) into $A_1 - B_1XB^{*}_1$ gives
\begin{align}
 A_1 - B_1XB^{*}_1 = A_1 - B_1X_0B^{*}_1  - B_1VF_BB^{*}_1 - B_1F_BV^{*}B^{*}_1 -
 B_1F_{B_2}U F_{B_3}B^{*}_1 - B_1F_{B_3}U^{*}F_{B_2}B^{*}_1,
 \label{33}
\end{align}
which is a matrix-valued function involving two variable matrices $V$ and $U$.
Thus, the constrained matrix-valued function in (\ref{yy11}) is equivalently
converted to the unconstrained matrix-valued function in (\ref{33}). To find
the global maximal and minimal ranks and partial inertias of (\ref{33}), we need the following result.

\begin{lemma} \label{T32}
Let
\begin{align}
 \phi(X_1, \, X_2) = A - B_1X_1C_1 - (B_1X_1C_1)^{*} - B_2X_2C_2 -
(B_2X_2C_2)^{*},
\label{34}
\end{align}
where $A\in \mathbb C_{{\rm H}}^{m},$ $B_i \in \mathbb C^{m \times p_i}$
and $C_i \in \mathbb C^{q_i \times m}$ are given$,$ and
 $X_i \in {\mathbb C}^{p_i \times q_i}$ are variable matrices for
 $i = 1, 2,$ and assume that
\begin{align}
{\mathscr R}(B_2) \subseteq {\mathscr R}(B_1), \ \ {\mathscr
R}(C^{*}_1) \subseteq {\mathscr R}(B_1), \ \ {\mathscr R}(C^{*}_2)
\subseteq {\mathscr R}(B_1).
\label{35}
\end{align}
 Also let
\begin{align*}
&  N =\left[\!\!\begin{array}{cccc} A & B_2 & C^{*}_1 & C^{*}_2 \\ C_1  & 0 & 0  & 0 \end{array}\!\!\right]\!, \ \
 N_1 =\left[\!\!\begin{array}{cccc} A & B_2 & C^{*}_1 & C^{*}_2 \\ B^{*}_2  & 0 & 0  & 0 \\ C_1  & 0 & 0  & 0
 \end{array}\!\!\right]\!,
\ \  N_2 =\left[\!\!\begin{array}{cccc} A & B_2 & C^{*}_1 & C^{*}_2 \\
C_1  & 0 & 0  & 0 \\ C_2 & 0 & 0  & 0
 \end{array}\!\!\right]\!,
\\
& M = \left[\!\! \begin{array}{cc} A  & B_1 \\ C_1 & 0 \end{array} \!\!\right], \ \
M_1 =\left[\!\!\begin{array}{ccc} A & B_2 & C^{*}_1  \\
 B_2^{*}  & 0 & 0 \\ C_1  & 0 & 0 \end{array}\!\!\right]\!, \ \ M_2 =\left[\!\!\begin{array}{ccc}
 A & C^{*}_1 & C^{*}_2  \\  C_1  & 0 & 0
\\ C_2  & 0 & 0   \end{array}\!\!\right]\!.
\end{align*}
Then$,$ the global maximal and minimal ranks and partial inertias of
$\phi(X_1, \, X_2)$ are given by
\begin{align}
\max_{X_1 \in {\mathbb C}^{p_1 \times q_1}, \, X_2 \in {\mathbb
C}^{p_2 \times q_2}} \!\!r[\, \phi(X_1, \, X_2)\,] & = \min\!\left\{
r[\, A,\, B_1\,], \ \
 r(N), \ \ r(M_1), \ \ r(M_2) \right\}\!,
 \label{36}
\\
\min_{X_1 \in {\mathbb C}^{p_1 \times q_1}, \, X_2 \in {\mathbb
C}^{p_2 \times q_2}} \!\!r[\, \phi(X_1, \, X_2)\,] & = 2r[\, A, \, B_1 \,] - 2r(M) + 2r(N) +
 \max\{\,  s_1, \ \  s_2, \ \  s_3, \ \  s_4 \,\},
\label{37}
\\
\max_{X_1 \in {\mathbb C}^{p_1 \times q_1}, \, X_2 \in {\mathbb
C}^{p_2 \times q_2}} \!\!i_{\pm}[\, \phi(X_1, \, X_2)\,] & =
\min\!\left\{ i_{\pm}(M_1), \ \ i_{\pm}(M_2) \right\}\!,
 \label{38}
\\
\min_{X_1 \in {\mathbb C}^{p_1 \times q_1}, \, X_2 \in {\mathbb
C}^{p_2 \times q_2}} \!\!i_{\pm}[\, \phi(X_1, \, X_2)\,] & = r[\, A,\, B_1\,]- r(M) + r (N) \nb
\\
&  \ \ \ \  + \max\{\,i_{\pm}(M_1)-r(N_1), \ \ i_{\pm}(M_2)-r(N_2)\, \},
  \label{39}
\end{align}
where
\begin{align*}
& s_1 = r(M_1) - 2r(N_1),  \ \ s_2 = r(M_2) - 2r(N_2),
\\
& s_3 = i_{+}(M_1) + i_{-}(M_2) - r(N_1) -  r(N_2),
\\
& s_4 = i_{+}(M_1) + i_{-}(M_2) - r(N_1) -  r(N_2).
\end{align*}
\end{lemma}

\noindent {\bf Proof} \ Under (\ref{35}), applying Lemma \ref{T18} to the variable matrix
$X_1$  in (\ref{34}) and simplifying, we
obtain
\begin{align}
\max_{X_1} r[\, \phi(X_1, \, X_2)\,] & = \min \left\{ r[\, A -
B_2X_2C_2 - (B_2X_2C_2)^{*}, \, B_1 \,], \ \
 r\!\left[\!\! \begin{array}{cc}  A - B_2X_2C_2 - (B_2X_2C_2)^{*}  & C_1^{*} \\ C_1 & 0 \end{array} \!\!\right]
 \right\} \nb
\\
& = \min \left\{ r[\, A, \, B_1 \,], \ \
 r\!\left[\!\! \begin{array}{cc}  A - B_2X_2C_2 - (B_2X_2C_2)^{*}  & C_1^{*} \\ C_1 & 0
 \end{array} \!\!\right] \right\},
\label{310}
\\
\min_{X_1} r[\, \phi(X_1, \, X_2)\,] & =  2r[\, A -  B_2X_2C_2 -
(B_2X_2C_2)^{*}, \, B_1 \, ] + r\! \left[\!\! \begin{array}{cc}  A -
B_2X_2C_2 - (B_2X_2C_2)^{*}  & C_1^{*} \\ C_1 & 0 \end{array}\!\!\right] \nb
\\
\max_{X_1} i_{\pm}[\, \phi(X_1, \, X_2)\,] & = i_{\pm}\!\left[\!\!
\begin{array}{cc}  A - B_2X_2C_2 - (B_2X_2C_2)^{*}  & C_1^{*} \\ C_1 & 0
\end{array} \!\!\right],
\label{311}
\\
&   \ \ \ - 2r\! \left[\!\! \begin{array}{cc}  A - B_2X_2C_2 -
(B_2X_2C_2)^{*} & B_1 \\ C_1 & 0 \end{array} \!\!\right] \nb
\\
&  =2r[\, A, \, B_1 \,] + r\! \left[\!\! \begin{array}{cc}  A -
B_2X_2C_2 - (B_2X_2C_2)^{*}
 & C_1^{*} \\ C_1 & 0 \end{array} \!\!\right]- 2r\! \left[\!\! \begin{array}{cc} A
  & B_1 \\ C_1 & 0 \end{array} \!\!\right],
\label{312}
\\
\min_{X_1}i_{\pm}[\, \phi(X_1, \, X_2)\,]  & = r[\, A -  B_2X_2C_2 -
(B_2X_2C_2)^{*}, \, B_1 \,] + i_{\pm}\!\left[\!\! \begin{array}{cc} A
- B_2X_2C_2 - (B_2X_2C_2)^{*}  & C_1^{*} \\ C_1 & 0 \end{array}
\!\!\right] \nb
\\
&   \ \ \ - r\! \left[\!\! \begin{array}{cc}  A - B_2X_2C_2 -
(B_2X_2C_2)^{*}  & B_1 \\ C_1 & 0 \end{array} \!\!\right] \nb
\\
&  =r[\, A, \, B_1 \,] + i_{\pm}\! \left[\!\! \begin{array}{cc}
A - B_2X_2C_2 - (B_2X_2C_2)^{*}
 & C_1^{*} \\ C_1 & 0 \end{array} \!\!\right]- r\!\left[\!\! \begin{array}{cc} A
  & B_1 \\ C_1 & 0 \end{array} \!\!\right].
\label{313}
\end{align}
Notice that
\begin{align}
\left[\!\!\begin{array}{cc} A - B_2X_2C_2 - (B_2X_2C_2)^{*}
 & C_1^{*} \\ C_1 & 0  \end{array} \!\!\right] & =
\left[\!\! \begin{array}{cc} A  &  C^{*}_1 \\ C_1 &  0
\end{array}\!\!\right]
 - \left[\!\!\begin{array}{c} B_2
\\  0 \end{array} \!\!\right]\!X_2[\, C_2, \, 0 \,] - \left[\!\! \begin{array}{c} C^{*}_2 \\  0 \end{array}
\!\!\right]X_2^{*}[\, B^{*}_2, \, 0 \,]　\nb
\\
& := \psi(X_2).
\label{314}
 \end{align}
Applying Lemma \ref{T15a} to this expression gives
\begin{align}
\max_{X_2\in \mathbb C^{m \times p_2}}\!\!r[\psi(X_2)] & = \min \left\{ \,
r(N), \ \ r(M_1), \ \ r(M_2) \, \right\}\!,
\label{316}
\\
\min_{X_2\in \mathbb C^{m \times p_2}}\!\!r[\psi(X_2)] & = 2r(N) + \max\{\,
 s_1, \ \  s_2, \ \ s_3, \ \ s_4 \, \},
\label{317}
\\
\max_{X_2\in \mathbb C^{m \times p_2}}\!\!i_{\pm}[\psi(X_2)] & =
\min\!\left\{ \, i_{\pm}(M_1), \ \ i_{\pm}(M_2)  \, \right\}\!,
\label{318}
\\
\min_{X_2\in \mathbb C^{m \times p_2}}\!\!i_{\pm}[\psi(X_2)]& = r(N) +
\max\{\, i_{\pm}(M_1) - r(N_1), \ \ i_{\pm}(M_2) - r(N_2) \, \},
\label{319}
\end{align}
where
$$
s_1 = r(M_1) - 2r(N_1),  \ \ s_2 = r(M_2) - 2r(N_2), \ \
$$
$$
s_3 = i_{+}(M_1) + i_{-}(M_2) - r(N_1) -  r(N_2), \ \
s_4 = i_{-}(M_1) + i_{+}(M_2) - r(N_1) -  r(N_2).
$$
Substituting these results into (\ref{310})--(\ref{313}) yields
(\ref{36})--(\ref{39}). \qquad $\Box$

\medskip

It is obviously of great importance to be able to give analytical
formulas for calculating the global maximal and minimal ranks and
inertias of the matrix expression in (\ref{34}) under the
assumptions in (\ref{36}). However, it is not easy to find the
global maximal and minimal ranks and inertias of a general
$\phi(X_1, \, X_2)$  as given in (\ref{34}). For convenience of
representation, we rewrite (\ref{33}) as
\begin{align}
A_1- B_1XB^{*}_1 = A - G_1VG_2 - (G_1VG_2)^{*} - G_3UG_4 - (G_3UG_4)^{*},
\label{319a}
\end{align}
where
\begin{align}
A = A_1 - B_1X_0B^{*}_1, \ \ G_1 =  B_1, \ \   G_2 = F_B B^{*}_1,  \ \
G_3 = B_1F_{B_2},  \ \ G_4 = F_{B_3}B^{*}_1.
\label{320}
\end{align}
It is easy to verify that the above matrices satisfy the conditions
\begin{align}
(G_2^{*}) \subseteq {\mathscr R}(G_1), \ \  {\mathscr R}(G_3) \subseteq {\mathscr R}(G_1),
\ \ {\mathscr R}(G_4^{*}) \subseteq {\mathscr R}(G_1), \ \ {\mathscr R}(G_2^{*}) \subseteq {\mathscr R}(G_3), \ \
{\mathscr R}(G_2^{*}) \subseteq {\mathscr R}(G^{*}_4).
\label{321}
\end{align}
In this case, applying Lemma \ref{T32} to (\ref{321}) yields the main results of
this section.

\begin{theorem}  \label{T33}
Let $A_i \in {\mathbb C}_{{\rm H}}^{m_i}$ and $B_i \in \mathbb
C^{m_i \times n}$ be given for $ i =1,\, 2,\, 3,$ and assume
 that the pair of matrix equations
 \begin{align}
 B_2 XB^{*}_2 = A_2 \ \  and  \ \ B_3 XB^{*}_3 = A_3
 \label{321a}
\end{align}
 have  a common solution $X \in {\mathbb C}_{{\rm H}}^{n}.$ Also denote the set of
 all their common Hermitian solutions by
\begin{align}
{\cal S} = \{ X \in {\mathbb C}_{{\rm H}}^{n} \ | \ B_2 XB^{*}_2 = A_2,  \ \  B_3 XB^{*}_3 = A_3 \, \}.
\label{321a1}
\end{align}
and let
\begin{align}
& P_1 =\left[\!\!\begin{array}{cccc} A_1 & B_1 & 0 & 0 \\ B_1^{*}  & 0 & B_2^{*}  & B_3^{*}
\end{array}\!\!\right], \ \ P_2 =\left[\!\!\begin{array}{ccc} A_1 &  0 &　B_1 \\　0 & -A_2 &
 B_2　\\  B^{*}_1 & B^{*}_2 & 0  \end{array}\!\!\right]\!, \ \
 P_3 =\left[\!\!\begin{array}{ccc} A_1 &  0 & B_1\\ 0 & -A_3 & B_3　
\\ B^{*}_1 & B^{*}_3 & 0  \end{array}\!\!\right]\!,
\label{322}
\\
& Q_1 =\left[\!\!\begin{array}{ccccc} A_1 & 0 & 0 & B_1 & B_1
  \\ 0 & -A_2 & 0  & B_2 & 0 \\ 0 & 0 & -A_3  & 0 & B_3 \\ B_1^{*} & B^{*}_2 & B^{*}_3 & 0 & 0
 \end{array}\!\!\right]\!, \  Q_2=\left[\!\!\begin{array}{cccc} A_1 &  0 &  B_1 & B_1 \\
0 &  -A_2 &  B_2 & 0 \\  B^{*}_1 &  B^{*}_2  & 0 & 0 \\ 0 & 0 & 0 & B_3
\end{array}\!\!\right]\!, \
 Q_3=\left[\!\!\begin{array}{cccc} A_1 &  0 &  B_1 & B_1 \\
0 & -A_3 &  B_3 & 0 \\  B^{*}_1 &  B^{*}_3  & 0 & 0 \\ 0 & 0 & 0 & B_2
\end{array}\!\!\right]\!.
\label{323}
\end{align}
  Then$,$
\begin{enumerate}
\item[{\rm(a)}] The global maximum rank of $A_1 - B_1XB^{*}_1$ subject to  {\rm (\ref{321a1})} is
\begin{align}
& \max_{X \in {\cal S}} r( \, A_1 - B_1XB^{*}_1 \,) \nb
\\
&  = \min\!\left\{r[\, A_1, \, B_1\,],  \ r(Q_1)- r \left[\!\!\begin{array}{c}
B_2  \\ B_3\end{array} \!\!\right] - r(B_2) - r(B_3),  \ r(P_2) -
2r(B_2), \  r(P_3) - 2r(B_3) \right\}.
 \label{324}
\end{align}

\item[{\rm(b)}] The global minimum rank of $A_1 - B_1XB^{*}_1$ subject to  {\rm (\ref{321a1})} is
\begin{align}
&  \min_{X \in {\cal S}} r( \, A_1 - B_1XB^{*}_1 \,) \nb
\\
&  = 2r[\, A_1, \, B_1 \,] - 2r(P_1) + 2r(Q_1) +
 \max\{ \, r(P_2) - 2r(Q_2), \  r(P_3) - 2r(Q_3), \  u_1, \  u_2 \,\},
\label{325}
\end{align}
where
$$
u_1 = i_{+}(P_2) + i_{-}(P_3) - r(Q_2) - r(Q_3), \ \ u_2 = i_{-}(P_2) +  i_{+}(P_3) - r(Q_2) - r(Q_3).
$$

\item[{\rm(c)}] The global maximum partial inertia of $A_1 - B_1XB^{*}_1$ subject to  {\rm (\ref{321a1})} is
\begin{align}
\max_{X \in {\cal S}} i_{\pm}( \, A_1 - B_1XB^{*}_1 \,) =\min \left\{ i_{\pm}(P_2) - r(B_2),
\ \ i_{\pm}(P_3) - r(B_3) \right\}.
\label{326}
\end{align}

\item[{\rm(d)}]  The global minimum partial inertia of $A_1 - B_1XB^{*}_1$ subject to  {\rm (\ref{321a1})} is
\begin{align}
 \min_{X \in {\cal S}} i_{\pm}(\, A_1 - B_1XB^{*}_1\,) & = r[\,A_1, \, B_1\,]- r(P_1) + r(Q_1) \nb
 \\
& \ \ \   + \max\!\left\{i_{\pm}(P_2) -  r(Q_2), \ i_{\pm}(P_3) -  r(Q_3) \right\}\!.
\label{327}
\end{align}
\end{enumerate}
\end{theorem}

\noindent {\bf Proof} \ Under (\ref{321}),  we find by Lemma \ref{T32} that
\begin{align}
 \max_{X \in {\cal S}} r(\,A_1 - B_1XB_1^{*} \,)
 & =\max_{V,\, U} r\left[ A - G_1VG_2 - (G_1VG_2)^{*}  - G_3UG_4 - (G_3UG_4)^{*} \right] \nb
  \\
& = \min \left\{ r[\, A, \, G_1 \,], \ r\!\left[\!\!\begin{array}{ccc}  A & G_3 & G_4^{*}
 \\ G_2 & 0 & 0 \end{array} \right]\!, \ r\!\left[\!\!\begin{array}{cc}
   A  & G_3  \\ G_3^{*} & 0  \end{array} \right]\!,  \
   r\!\left[\!\!\begin{array}{cc} A & G^{*}_4 \\ G_4 & 0 \end{array} \!\!\right] \right\}\!,
\label{328}
\\
 \min_{X \in {\cal S}} r(\,A_1 - B_1XB_1^{*} \,)
& =\min_{V,\, U} r\left[\,A - G_1VG_2 - (G_1VG_2)^{*}  - G_3UG_4 - (G_3UG_4)^{*} \, \right] \nb
  \\
& = 2r[\, A, \, G_1 \,] - 2r\left[\!\!\begin{array}{cc}
   A & G_1  \\  G_2 & 0 \end{array} \!\!\right] + 2r\left[\!\!\begin{array}{ccc}
   A & G_3  & G_4^{*}
   \\ G_2 & 0 & 0\end{array} \!\!\right]  \nb
   \\
   & \ \ \ \ + \max\{ s_1, \ \ s_2, \ \  s_3, \ \ s_4 \, \},
\label{329}
\\
 \max_{X \in {\cal S}} i_{\pm}(\,A_1 - B_1XB_1^{*} \,) & =
\max_{V,\, U}  i_{\pm}\left[ A - G_1VG_2 - (G_1VG_2)^{*}  - G_3UG_4 - (G_3UG_4)^{*} \right] \nb
  \\
& = \min \left\{  i_{\pm}\!\left[\!\!\begin{array}{cc}
   A  & G_3  \\ G_3^{*} & 0  \end{array} \right]\!,  \
   i_{\pm}\!\left[\!\!\begin{array}{cc} A & G^{*}_4 \\ G_4 & 0 \end{array}\!\!\right] \right\}\!,
\label{330}
\\
 \min_{X \in {\cal S}} i_{\pm}(\,A_1 - B_1XB_1^{*} \,) & =\min_{V,\, U}  i_{\pm}\left[\,A - G_1VG_2 - (G_1VG_2)^{*}
 - G_3UG_4 - (G_3UG_4)^{*} \, \right] \nb
  \\
& = r[\, A, \, G_1 \,] - r\left[\!\!\begin{array}{cc}
   A & G_1  \\  G_2 & 0 \end{array} \!\!\right] + r\left[\!\!\begin{array}{ccc}
   A & G_3  & G_4^{*}
   \\ G_2 & 0 & 0\end{array} \!\!\right] + \max\{ t_1, \ \ t_2 \, \},
\label{331}
\end{align}
where
\begin{align*}
s_1 & = r\!\left[\!\!\begin{array}{cc}  A & G_3  \\ G_3^{*} & 0 \end{array} \!\!\right] -
   2r\!\left[\!\!\begin{array}{cccc}  A & G_3 & G_4^{*}
    \\
   G_3^{*} & 0 & 0 \end{array} \!\!\right]\!,
   \\
s_2 & = r\!\left[\!\!\begin{array}{ccc}  A & G_4^{*} \\ G_4 & 0
\end{array}\!\! \right] - 2r\!\left[\!\!\begin{array}{cccc}
   A & G_3 & G_4^{*}  \\ G_4 & 0 & 0 \end{array} \!\!\right]\!,
\\
s_3 &  = i_ {+}\!\left[\!\!\begin{array}{ccc}  A & G_3  \\ G_3^{*} & 0 \end{array} \right]
 + i_{-}\!\left[\!\!\begin{array}{ccc}
   A & G_4^{*}  \\   G_4 & 0  \end{array} \!\!\right] -
    r\!\left[\!\!\begin{array}{ccc}  A & G_3 & G_4^{*}  \\  G_3^{*} & 0 & 0 \end{array}\!\! \right]
    - r\left[\!\!\begin{array}{cccc}
   A & G_3 & G_4^{*}  \\
 G_4 & 0 & 0
\end{array}\!\! \right]\!,
\\
s_4 & =  i_{-}\!\left[\!\!\begin{array}{ccc}
   A  & G_3  \\  G_3^{*} & 0  \end{array} \!\!\right]
    + i_{+}\left[\!\!\begin{array}{cc}  A &  G_4^{*}  \\
   G_4 & 0\end{array}\right] - r\!\left[\!\!\begin{array}{ccc}
   A  & G_3 & G_4^{*}  \\  G_3^{*} & 0 & 0 \end{array}\!\!\right]
    - r\left[\!\begin{array}{ccc}  A  & G_3 & G_4^{*}  \\
   G_4 & 0 & 0 \end{array} \!\!\right]\!,
\\
t_1 & = i_{\pm}\!\left[\!\!\begin{array}{cc}  A & G_3 \\
   G_3^{*} & 0 \end{array} \!\!\right] - r\!\left[\!\!\begin{array}{ccc}
   A & G_3 & G_4^{*}  \\ G_3^{*} & 0 & 0 \end{array} \!\!\right]\!,
   \\
t_2 & = i_{\pm}\!\left[\!\!\begin{array}{cc}  A & G_4 \\
   G_4^{*} & 0 \end{array} \right] - r\!\left[\!\!\begin{array}{ccc}
   A & G_3 & G_4^{*} \\ G_4^{*} & 0 & 0 \end{array}\!\!\right]\!.
   \end{align*}
Applying  (\ref{17})--(\ref{19}) and (\ref{114u}), and simplifying by
$[\,B_2X_0B^{*}_2, \, B_3X_0B_3^{*}\,] = [\,A_2, \, A_3\,]$,  elementary matrix operations and
 congruence matrix operations,  we obtain
\begin{align}
 r[\, A, \, G_1 \,] & = r[\, A_1 - B_1X_0B_1^{*}, \, B_1 \,] = r[\, A_1, \, B_1 \,],
\label{332}
\\
r\!\left[\!\!\begin{array}{ccc}  A & G_3 & G_4^{*} \\ G_2 & 0 & 0 \end{array}\!\! \right]
&  = r\!\left[\!\!\begin{array}{ccc}
   A_1 - B_1X_0B_1^{*}  & B_1F_{B_2} & B_1 F_{B_3}
   \\
   F_B B_1^{*} & 0 & 0 \end{array}\!\!\right] \nb
\\
& = r\left[\!\!\begin{array}{cccc}
   A_1 - B_1X_0B_1^{*} & B_1 & B_1 & 0
   \\
   B_1^{*} & 0 & 0 & B^{*}
   \\
   0 & B_2 & 0 & 0
   \\
 0 & 0 & B_3 & 0\end{array}\!\!\right] - r(B) - r(B_2) - r(B_3) \nb
 \\
& = r\left[\!\!\begin{array}{cccc} A_1 & B_1 & B_1 &  B_1X_0B^{*}
   \\
   B_1^{*} & 0 & 0 & B^{*}
   \\
   0 & B_2 & 0 & 0
   \\
 0 & 0 & B_3 & 0\end{array}\!\!\right] - r(B) - r(B_2) - r(B_3) \nb
 \\
& = r\left[\!\!\begin{array}{ccccc} A_1 & B_1 & B_1 &  0 & 0
   \\
   B_1^{*} & 0 & 0 & B^{*}_2 & B^{*}_3
   \\
   0 & B_2 & 0 &  - A_2 & 0
   \\
 0 & 0 & B_3 & 0 &  - A_3\end{array}\!\!\right] - r(B) - r(B_2) - r(B_3) \nb
 \\
&  = r(Q_1) - r(B) - r(B_2) - r(B_3),
\label{333}
\\
r\!\left[\!\!\begin{array}{cc}
   A & G_1   \\  G_2 & 0 \end{array} \!\!\right] & = r\!\left[\!\!\begin{array}{cc}
   A_1 - B_1X_0B_1^{*} & B_1
    \\
   F_B B_1^{*} & 0 \end{array}\!\! \right]  = r\!\left[\!\!\begin{array}{ccc} A_1  & B_1 & 0
    \\
    B_1^{*} & 0  & B^{*} \end{array}\!\! \right]  - r(B) \nb
\\
  & = r(P_1) - r(B),
\label{334}
\\
i_{\pm}\!\left[\!\!\begin{array}{cc}  A & G_3  \\ G^{*}_3 & 0 \end{array}\!\! \right]
&  = i_{\pm}\!\left[\!\!\begin{array}{cc}
   A_1 - B_1X_0B_1^{*}  & B_1F_{B_2}
   \\
   F_{B_2}B_1^{*} & 0 \end{array}\!\!\right] = i_{\pm}\!\left[\!\!\begin{array}{ccc}
   A_1 - B_1X_0B_1^{*} & B_1 & 0
    \\
   B^{*}_1  & 0 & B^{*}_2
   \\
  0 & B_2 & 0\end{array}\!\!\right]- r(B_2) \nb
\\
& = i_{\pm}\!\left[\!\!\begin{array}{ccc}
   A_1 & B_1 &  B_1X_0B_2^{*}/2
    \\
   B^{*}_1  & 0 & B^{*}_2
   \\
  B_1X_0B_2^{*}/2 & B_2 & 0\end{array}\!\!\right]- r(B_2)  = i_{\pm}\!\left[\!\!\begin{array}{ccc}
   A_1  & B_1 & 0
    \\
   B^{*}_1  & 0 & B^{*}_2
   \\
  0 & B_2 & -A_2\end{array}\!\!\right]- r(B_2) \nb
  \\
  & = i_{\pm}(P_2) - r(B_2),
\label{335}
\end{align}
\begin{align}
 r\!\left[\!\!\begin{array}{ccc}
   A & G_3 & G_4^{*}
   \\
   G^{*}_3 & 0 & 0
\end{array}\!\!\right] & = r\!\left[\!\!\begin{array}{ccc}
   A_1 - B_1X_0B_1^{*} & B_1 F_{B_2} & B_1 F_{B_3}
   \\
   F_{B_2}B_1^{*} & 0 & 0
\end{array}\!\!\right] \nb
 \\
 & = r\!\left[\!\!\begin{array}{ccccc}
   A_1 - B_1 X_0B_1^{*} & B_1 & B_1 & 0
   \\
   B_1^{*} & 0 & 0 & B_2^{*}
   \\
    0 & B_2 & 0 & 0
   \\
   0 & 0 & B_3 & 0
\end{array} \!\!\right] - 2r(B_2) - r(B_3) \nb
\\
 & = r\!\left[\!\!\begin{array}{cccc}
   A_1  & B_1 & B_1 & B_1 X_0B^{*}_2
   \\
   B_1^{*} & 0 & 0 & B_2^{*}
   \\
    0 & B_2 & 0 & 0
   \\
   0 & 0 & B_3 & 0
\end{array} \!\!\right] - 2r(B_2) - r(B_3) \nb
\\
& = r\!\left[\!\!\begin{array}{ccccc}
   A_1  & B_1 & B_1 & 0
   \\
   B_1^{*} & 0 & 0 & B_2^{*}
   \\
    0 & B_2 & 0 & -A_2
   \\
   0 & 0 & B_3 & 0
\end{array} \!\!\right] - 2r(B_2) - r(B_3) \nb
\\
&  = r(Q_2) - 2r(B_2) - r(B_3).
 \label{336}
\end{align}
By a similar approach, we can obtain
\begin{align}
i_{\pm}\!\left[\!\!\begin{array}{cc}  A & G_4  \\ G^{*}_4 & 0 \end{array}\!\! \right]  = i_{\pm}(P_3)
- r(B_3), \ \ \
 r\!\left[\!\!\begin{array}{ccc}
   A & G_3 & G_4^{*}
   \\
   G_4 & 0 & 0
\end{array}\!\!\right] = r(Q_3) - r(B_2) - 2r(B_3).
\label{338}
\end{align}
Substituting (\ref{332})--(\ref{338}) into (\ref{328})--(\ref{331}) yields
(\ref{324})--(\ref{327}). \qquad $\Box$

\medskip

Some direct consequences of the previous theorem are given below.

\begin{corollary}  \label{T34}
Let $A_i \in \mathbb C_{{\rm H}}^{m_i}$ and $B_i \in \mathbb C^{m_i \times n}$ be given
for $ i =1,\, 2,\, 3,$ and suppose that each pair of $B_1 XB^{*}_1 = A_1,$ $B_2 XB^{*}_2 = A_2$
and $B_3 XB^{*}_3 = A_3$ have a common Hermitian solution$.$  Also let  ${\cal S}$ be
of the form {\rm (\ref{321a})}$.$   Then$,$
\begin{align}
\max_{X \in {\cal S}} r( \, A_1 - B_1XB^{*}_1 \,) & = \min\!\left\{r(B_1),  \ r(Q_1) - r\!\left[\!\!\begin{array}{c}
B_2  \\ B_3\end{array} \!\!\right] - r(B_2) - r(B_3), \right. \nb
\\
& \ \ \ \ \ \ \ \ \ \ \ \left. 2r\!\left[\!\!\begin{array}{c}
B_1  \\ B_2\end{array} \!\!\right]  - 2r(B_2), \ 2r\!\left[\!\!\begin{array}{c}
B_1  \\ B_3\end{array} \!\!\right]  - 2r(B_3) \right\},
 \label{340}
\\
 \min_{X \in {\cal S}} r( \, A_1 - B_1XB^{*}_1 \,) & = 2r(Q_1) -  2r\!\left[\!\!\begin{array}{c}
B_1  \\ B_2 \\ B_3 \end{array} \!\!\right] - 2r\!\left[\!\!\begin{array}{cccc} B_1 & B_1 \\
B_2 & 0 \\ 0 & B_3 \end{array}\!\!\right]\!,
\label{341}
\\
 \max_{X \in {\cal S}} i_{\pm}( \, A_1 - B_1XB^{*}_1 \,) & =\min \left\{ r\!\left[\!\!\begin{array}{c}
B_1  \\ B_2\end{array} \!\!\right] - r(B_2), \ \ r\!\left[\!\!\begin{array}{c}
B_1  \\ B_3\end{array} \!\!\right] - r(B_3) \right\},
\label{342}
\\
 \min_{X \in {\cal S}} i_{\pm}( \, A_1 - B_1XB^{*}_1 \,) & = r(Q_1) - r\!\left[\!\!\begin{array}{c}
B_1  \\ B_2 \\ B_3 \end{array} \!\!\right] - r\!\left[\!\!\begin{array}{cccc} B_1 & B_1 \\
B_2 & 0 \\ 0 & B_3 \end{array}\!\!\right]\!,
 \label{343}
\end{align}
where $Q_1$ is of the form {\rm (\ref{323})}$.$

\end{corollary}

\noindent {\bf Proof} \  Under the given conditions, the ranks and inertias of the block
matrices in (\ref{322}) and (\ref{323}) are given by
$$
r(P_1) = r(B_1) + r\!\left[\!\!\begin{array}{c}
B_1  \\ B_2 \\ B_3\end{array} \!\!\right], \ \ r(P_2) = 2r\!\left[\!\!\begin{array}{c}
B_1  \\ B_2\end{array} \!\!\right], \
r(P_3) = 2r\!\left[\!\!\begin{array}{c}
B_1  \\ B_3\end{array} \!\!\right], \
 i_{\pm}(P_2) = r\!\left[\!\!\begin{array}{c}
B_1  \\ B_2\end{array} \!\!\right], \  i_{\pm}(P_3) = r\!\left[\!\!\begin{array}{c}
B_1  \\ B_3\end{array} \!\!\right],
$$
$$
r(Q_2) = r\left[\!\!\begin{array}{cc} B_1 & B_1 \\
  B_2 & 0 \\  0 & B_3
\end{array}\!\!\right] + r\!\left[\!\!\begin{array}{c}
B_1  \\ B_2\end{array} \!\!\right], \ \  r(Q_3) = r\left[\!\!\begin{array}{cc} B_1 & B_1 \\
  B_2 & 0 \\ 0 & B_3
\end{array}\!\!\right] + r\!\left[\!\!\begin{array}{c}
B_1  \\ B_3\end{array} \!\!\right]\!.
$$
Hence (\ref{324})--(\ref{327}) reduce to (\ref{340})--(\ref{343}). \qquad $\Box$

\begin{corollary}  \label{T35}
Let $A_i \in \mathbb C_{{\rm H}}^{m_i \times m_i}$ and $B_i \in \mathbb
C^{m_i \times n}$ be given for $ i =1,\, 2,\, 3,$ and suppose
 that each pair of the triple matrix equations
\begin{align}
B_1 XB^{*}_1 = A_1, \ \ B_2 XB^{*}_2 = A_2, \ \ B_3 XB^{*}_3 = A_3
 \label{345}
\end{align}
have a common Hermitian solution$.$ Then$,$ there exists a Hermitian $X$ such that {\rm (\ref{345})}
holds  if and only if
\begin{align}
r\!\left[\!\!\begin{array}{ccccc} A_1 & 0 & 0 & B_1 & B_1
  \\ 0 & -A_2 & 0  & B_2 & 0 \\ 0 & 0 & -A_3  & 0 & B_3 \\ B_1^{*} & B^{*}_2 & B^{*}_3 & 0 & 0
 \end{array}\!\!\right] = r\!\left[\!\!\begin{array}{cccc} B_1 & B_1 \\
B_2 & 0 \\ 0 & B_3 \end{array}\!\!\right]  + r[\, B^{*}_1,  \, B^{*}_2, \, B^{*}_3 \,].
 \label{346}
\end{align}
\end{corollary}

\noindent {\bf Proof} \ It follows from (\ref{341}). \qquad $\Box$

\medskip

A challenging open problem on the triple matrix equations in (\ref{345}) is to give a parametric form for
their general common Hermitian solution.

Setting $B_1 = I_n$ in Theorem 3.3 may yield a group of results on the extremal
ranks and inertias of $A_1 - X$ subject to {\rm (\ref{321a1})}. In particular, we have
the following consequences.

\begin{corollary}  \label{T36}
Let $A_i \in \mathbb C_{{\rm H}}^{m_i}$ and $B_i \in \mathbb
C^{m_i \times n}$ be given for $i = 2,\, 3,$ and assume that
{\rm (\ref{321a})} has a common solution$.$  Also let ${\cal S}$ be of the form
{\rm (\ref{321a1})}$.$ Then$,$
\begin{enumerate}
\item[{\rm (a)}] The global maximum rank of the solution of {\rm (\ref{321a1})}
is
\begin{align}
 \max_{X \in {\cal S}}  r(X)  = \min\{\, n, \ \ s_1, \ \ s_2, \ \ s_3 \, \},
 \label{347}
\end{align}
where
$$
s_1 = 2n + r\!\left[\!\!\begin{array}{ccc} A_2 & 0 & B_2 \\
0 & A_3 & B_3 \end{array}\!\!\right]  - r\!\left[\!\!\begin{array}{c}
B_2  \\ B_3\end{array} \!\!\right] - r(B_2) - r(B_3),
$$
$$
s_2= 2n + r(A_2) - 2r(B_2), \ \ s_3 = 2n + r(A_3) - 2r(B_3).
$$

\item[{\rm (b)}] The global minimum rank of the solution of {\rm (\ref{321a1})} is
\begin{align}
 \min_{X \in {\cal S}} r(X) & =
 2r\!\left[\!\!\begin{array}{ccc} A_2 & 0 & B_2 \\
0 & A_3 & B_3 \end{array}\!\!\right]  + \max\{ \, t_1, \ \ t_2, \ \  t_3, \ \ t_4 \, \},
 \label{348}
\end{align}
where
\begin{align*}
t_1 & =  r(A_2) - 2r\!\left[\!\!\begin{array}{ccc} A_2 & B_2 \\
0 & B_3 \end{array}\!\!\right]\!, \ \ \ t_2 = r(A_3) - 2r\!\left[\!\!\begin{array}{ccc} 0 & B_2 \\
A_3 & B_3 \end{array}\!\!\right]\!,  \ \
\\
t_3 & = i_{+}(A_2) + i_{-}(A_3) - r\!\left[\!\!\begin{array}{ccc} A_2 & B_2 \\
0 & B_3 \end{array}\!\!\right] - r\!\left[\!\!\begin{array}{ccc} 0 & B_2 \\
A_3 & B_3 \end{array}\!\!\right]\!, \ \
\\
t_4 & = i_{-}(A_2) + i_{+}(A_3) - r\!\left[\!\!\begin{array}{ccc} A_2 & B_2 \\
0 & B_3 \end{array}\!\!\right] - r\!\left[\!\!\begin{array}{ccc} 0 & B_2 \\
A_3 & B_3 \end{array}\!\!\right]\!.
\end{align*}

\item[{\rm (c)}] The global maximum partial inertia of the solution of {\rm (\ref{321a1})}  is
\begin{align}
 \max_{X \in {\cal S}}  i_{\pm}(X)
  & = \min\{ \, n +  i_{\pm}(A_2) - r(B_2), \ \  n +  i_{\pm}(A_3) -  r(B_3)\,
  \}.
\label{349}
\end{align}

\item[{\rm (d)}] The global minimum partial inertia of the solution of
{\rm (\ref{321a1})} is
\begin{align}
 \min_{X \in {\cal S}}  i_{\pm}
  (X) & = r\!\left[\!\!\begin{array}{ccc} A_2 & 0 & B_2 \\
0 & A_3 & B_3 \end{array}\!\!\right]  + \max \left\{ i_{\pm}(A_2) -
r\!\left[\!\!\begin{array}{ccc} A_2 & B_2 \\
0 & B_3 \end{array}\!\!\right]\!, \ i_{\pm}(A_3) -
r\!\left[\!\!\begin{array}{ccc} 0 & B_2 \\
A_3 & B_3 \end{array}\!\!\right] \! \right\}\!.
 \label{350}
\end{align}
\end{enumerate}
In consequence$,$
\begin{enumerate}
\item[{\rm (e)}]  Eq.\,{\rm (\ref{321a})} has a solution  $X >0$ if
and only if
$$
A_2 \geqslant 0, \ \ A_3 \geqslant 0, \  {\mathscr R}(A_2) = {\mathscr R}(B_2),  \ \
{\mathscr R}(A_3) = {\mathscr R}(B_3).
$$

\item[{\rm (f)}] All solutions of {\rm (\ref{321a})} satisfy $X >0$ if
and only if  $A_2 \geqslant 0$, $A_3 \geqslant 0$ and one of
$$
 r(A_2) = r(B_2) = n, \ \ \  r(A_3) = r(B_3) =n.
$$

\item[{\rm (g)}]  Eq. {\rm (\ref{321a})} has a solution  $X <0$ if
and only if
$$
A_2 \leqslant 0, \ \ A_3 \leqslant 0, \ \ {\mathscr R}(A_2) = {\mathscr R}(B_2),  \ \ {\mathscr R}(A_3)
=  {\mathscr R}(B_3).
$$

\item[{\rm (h)}] All solutions of {\rm (\ref{321a})}  satisfy $X <0$ if
and only if $A_2 \leqslant 0$, $A_3 \leqslant 0$ and one of
$$
 r(A_2) = r(B_2) = n, \ \ \  r(A_3) = r(B_3) =n.
$$

\item[{\rm (i)}]  Eq. {\rm (\ref{321a})} has a solution $X \geqslant 0$ if
and only if
$$
A_2 \geqslant 0, \ \ A_3 \geqslant 0, \ \ {\mathscr R}\!\left[\!\!\begin{array}{c} A_2 \\
0  \end{array}\!\!\right] \subseteq  {\mathscr R}\!\left[\!\!\begin{array}{cc} 0 & B_2 \\
A_3 & B_3 \end{array}\!\!\right]\!, \ \
{\mathscr R}\!\left[\!\!\begin{array}{c} 0 \\
A_3  \end{array}\!\!\right] \subseteq  {\mathscr R}\!\left[\!\!\begin{array}{ccc} A_2 & B_2 \\
0 & B_3 \end{array}\!\!\right]\!.
$$

\item[{\rm (j)}] All solutions of {\rm (\ref{321a})} satisfy $X \geqslant 0$ if
and only if $A_2 \geqslant 0,$ $A_3 \geqslant 0$ and one of
$$
r(B_2) =  n \ \ and \ \  r(B_3) =  n.
$$

\item[{\rm (k)}]  Eq. {\rm (\ref{321a})} has a solution  $X \leqslant 0$ if
and only if
$$
A_2 \leqslant 0, \ \ A_3 \leqslant 0, \ \ {\mathscr R}\!\left[\!\!\begin{array}{c} A_2 \\
0  \end{array}\!\!\right] \subseteq  {\mathscr R}\!\left[\!\!\begin{array}{cc} 0 & B_2 \\
A_3 & B_3 \end{array}\!\!\right]\!, \ \
{\mathscr R}\!\left[\!\!\begin{array}{c} 0 \\
A_3  \end{array}\!\!\right] \subseteq  {\mathscr R}\!\left[\!\!\begin{array}{ccc} A_2 & B_2 \\
0 & B_3 \end{array}\!\!\right]\!.
$$

\item[{\rm (l)}] All solutions of {\rm (\ref{321a})} satisfy $X \leqslant 0$ if
and only if $A_2 \leqslant 0$, $A_3 \leqslant 0$ and one of
$$
r(B_2) =  n \ \ and \ \  r(B_3) =  n.
$$
\end{enumerate}
\end{corollary}

\noindent {\bf Proof} \ Set $A_1 = 0$ and $B_1 = I_n$ in Theorem \ref{T33} and simplifying,
we obtain (a)--(d). Applying Lemma \ref{T12} to  (\ref{348}) and (\ref{349}), we obtain
(e)--(l). \qquad $\Box$

\medskip

Corollary \ref{T36}(e)--(l) give a set of analytical characterizations for the existence of definite
common solutions of the two matrix equations in (\ref{321a}) by using some rank and range equalities
and inequalities. These characterizations are simple and easy to understand in comparison with
various known conditions (see, e.g., \cite{LR,Zh1,Zh2})s on the existence of definite common solutions
of (\ref{321a}).

Rewrite  $B_2 XB^{*}_2 = A_2$ and $B_3 XB^{*}_3 = A_3$ as
\begin{equation}
  [\, B_{21},  \,  B_{22} \,]\! \left[\!\! \begin{array}{cc}  X_1  & X_2
\\ X_2^{*} & X_3 \end{array} \!\!\right]\!\left[\!\! \begin{array}{c} B^{*}_{21} \\ B^{*}_{22}
\end{array} \!\!\right] = A_2, \ \
 [\, B_{31},  \,  B_{32} \,]\! \left[\!\! \begin{array}{cc}  X_1  & X_2
\\ X_2^{*} & X_3 \end{array} \!\!\right]\!\left[\!\! \begin{array}{c} B^{*}_{31} \\ B^{*}_{32}
\end{array} \!\!\right] = A_3,
\label{352}
\end{equation}
where $B_{i1} \in {\mathbb C}^{m_i \times n_1},$ $B_{i2} \in \mathbb
C^{m_i \times n_2},$ $i =2, \, 3,$ $X_1 \in {\mathbb C}_{{\rm H}}^{n_1}, \ X_2 \in
\mathbb C^{n_1 \times n_2}$ and $X_3 \in  \mathbb C_{{\rm H}}^{n_2}$
 with $n_1 + n_2 = n$.  We next derive the extremal ranks and inertias
 of the submatrices $X_1$ and $X_3$ in a Hermitian solution of (\ref{352}).
 Note that $X_1, \, X_2, \, X_3$ in (\ref{352}) can be rewritten as
\begin{align}
 X_1 = P_1XP^{*}_1, \ \ X_2 = P_1XP^{*}_2, \ \ X_3= P_2XP^{*}_2,
\label{353}
\end{align}
where $P_1 = [\, I_{n_1},  \,  0 \,]$ and $P_2 =[ \, 0, \, I_{n_2} \, ].$
For convenience, we adopt the following notation for the collections of the
submatrices $X_1$ and $X_3$  in (\ref{352}):
\begin{align}
 {\cal S}_1 & =  \left\{ X_1 = P_1XP^{*}_1 \ | \ B_2 XB^{*}_2 = A_2, \ B_3 XB^{*}_3 = A_3, \ X = X^{*}  \right\},
\label{354}
\\
 {\cal S}_3 & = \left\{ X_3 = P_2XP^{*}_2 \ | \ B_2 XB^{*}_2 = A_2, \ B_3 XB^{*}_3 = A_3, \ X = X^{*}  \right\}.
\label{355}
\end{align}
The global maximal and minimal ranks and partial inertias of the submatrices $X_1$ and $X_3$  in (\ref{352})
can easily be derived from Theorem \ref{T33}. The details are omitted.

If each of the triple matrix equations in (\ref{yy17}) is not consistent,
people may alternatively seek its common approximation solutions under
various given optimal criteria. One of the most useful approximation
solutions of $BXB^{*} = A$ is the well-known least-squares Hermitian solution,
which is defined to be a Hermitian matrix $X$ that minimizes the objective function:
\begin{equation}
\|\,A - BXB^{*}\,\|^2 = {\rm tr}[\, (\,A - BXB^{*}\,)(\,A - BXB^{*}\,)^{*} \,].
\label{ff355}
\end{equation}
The normal equation corresponding to the norm minimization problem is given by
\begin{equation}
B^{*}BXB^{*}B = B^{*}AB.
\label{ff356}
\end{equation}
This equation is always consistent. Concerning the common least-squares
Hermitian solution of (\ref{yy17}), we have the following result.

\begin{corollary}  \label{T37}
Let $A_i \in \mathbb C_{{\rm H}}^{m_i}$ and $B_i \in \mathbb
C^{m_i \times n}$ be given for $i =1,\, 2,\, 3.$ Then$,$ triple matrix
equations have a common least-squares Hermitian solution$,$ namely$,$  there
exists an $X \in \mathbb C_{{\rm H}}^{n \times n}$ such that
 \begin{equation}
\|\,A_i - B_iXB^{*}_i\,\| = \min,  \ \ i =1,\, 2,\, 3,
\label{357}
\end{equation}
if and only if
\begin{align}
&r\!\left[\!\!\begin{array}{ccc} B_i^{*}A_iB_i &  0 &　B_i^{*}B_i \\　0 & -B_j^{*}A_jB_j &
 B_j^{*}B_j　\\  B^{*}_iB_i & B^{*}_jB_j & 0  \end{array}\!\!\right] = 2r\!\left[\!\!\begin{array}{c}
B_i  \\ B_j \end{array} \!\!\right]\!, \ \ i \neq j, \  i, \, j = 1,\, 2,\, 3,
 \label{ff358}
 \\
& r\!\left[\!\!\begin{array}{ccccc} B_1^{*}A_1B_1 & 0 & 0 & B_1^{*}B_1 & B_1^{*}B_1
  \\ 0 & -B_2^{*}A_2B_2 & 0  & B_2^{*}B_2 & 0 \\ 0 & 0 & -B_3^{*}A_3B_3  & 0 & B_3^{*}B_3 \\
   B_1^{*}B_1 & B^{*}_2B_2 & B^{*}_3B_3 & 0 & 0
 \end{array}\!\!\right] = r\!\left[\!\!\begin{array}{cccc} B_1 & B_1 \\
B_2 & 0 \\ 0 & B_3 \end{array}\!\!\right]  + r\!\left[\!\!\begin{array}{c}
B_1  \\ B_2 \\ B_3 \end{array} \!\!\right]\!.
 \label{ff359}
\end{align}
\end{corollary}

\noindent {\bf Proof} \ It follows from Lemma \ref{T31}, Corollary \ref{T35} and
(\ref{ff356}). \qquad $\Box$

\section{The extremal ranks and inertias of  $A_1-B_1XB_1^{*}$ subject to  the Hermitian solutions of $B_4X = A_4$}
\setcounter{section}{4} \setcounter{equation}{0}

Also $B_4X = A_4$ in (\ref{yy12}) is not given in symmetric pattern, it may have a Hermitian solution, as shown
in Theorem \ref{T15b}. So that the global extremal ranks and inertias of $A_1-B_1XB_1^{*}$ subject to
the Hermitian solution or nonnegative definite solution of $B_4X = A_4$ can also be derived.

\begin{theorem} \label{T41}
Assume that the matrix equation $B_4X = A_4$ in {\rm (\ref{yy12})}
has a Hermitian solution$,$  i.e.$,$ ${\mathscr R}(A_4)\subseteq
{\mathscr R}(B_4)$ and $A_4B^{*}_4 =B_4A^{*}_4,$ and let
\begin{align}
{\mathcal S} = \{\, X\in {\mathbb C}^{n}_{{\rm H}} \ | \  B_4X = A_4 \, \}, \ \
M = \left[\! \begin{array}{cc}  A_1 & B_1 \\
A_4B_1^{*} & B_4 \end{array} \!\right], \ \  N =\left[\! \begin{array}{ccc}  A_1 & B_1 & 0
\\
B_1^{*} & 0 & B_4^{*} \\ 0 & B_4 &  -A_4B_4^{*}  \end{array} \!\right].
\label{41}
\end{align}
Then$,$
 \begin{align}
\max_{X\in {\mathcal S}} r(\,A_1-B_1XB_1^{*}\,) & = r(M) - r(B_4),
 \label{42}
\\
\min_{X\in {\mathcal S}} r(\,A_1-B_1XB_1^{*}\,) & = 2r(M)- r(N),
 \label{43}
\\
\max_{X\in {\mathcal S}}
i_{\pm}(\,A_1-B_1XB_1^{*}\,) & = i_{\pm}(N)  - r(B_4),
 \label{44}
\\
\min_{X\in {\mathcal S}}
i_{\pm}(\,A_1-B_1XB_1^{*}\,) & = r(M) -  i_{\mp}(N).
 \label{45}
\end{align}
In  consequences$,$
\begin{enumerate}
\item[{\rm (a)}] $B_4X=A_4$ has a solution $X \in \mathbb C^n_{{\rm H}}$ such that $A_1-B_1XB_1^{*}$  is nonsingular
 if and only if $r(M)=r(B_4)+m_1.$

\item[{\rm (b)}]  $A_1-B_1XB_1^{*}$  is nonsingular for all Hermitian solution of $B_4X=A_4$
 if and only if $2r(M)=r(N) + m_1.$

\item[{\rm (c)}] The pair of matrix equations  $B_1XB_1^{*} = A_1$ and $B_4X=A_4$ have a common Hermitian solution
 if and only if ${\mathscr R}\!\left[\!\! \begin{array}{cc}  A_1   \\ A_4B_1^{*} \end{array} \!\!\right] \subseteq
{\mathscr R}\!\left[\!\! \begin{array}{cc}   B_1  \\ B_4 \end{array}
\!\!\right].$

\item[{\rm (d)}]  $B_1XB_1^{*} = A_1$ holds for all Hermitian solutions of $B_4X=A_4$
 if and only if $r(M) = r(B_4).$

\item[{\rm (e)}] $B_4X=A_4$ has a solution $X \in {\mathbb C}^n_{\rm H}$ such that $ A_1 - B_1XB_1^{*} >0$
 $(\, A_1 - B_1XB_1^{*} <0 \,)$ if and only if $i_{+}(N) = r(B_4) + m_1 \ \ (\, i_{-}(N)= r(B_4) + m_1 \,).$

\item[{\rm (f)}]  $A_1 - B_1XB_1^{*} >0$ $(\, A_1 - B_1XB_1^{*} <0\,)$  holds for all Hermitian solutions of $B_4X=A_4$
 if and only if $r(M)= i_{-}(N) + m_1$ $(\, r(M)= i_{+}(N) + m_1\,).$

\item[{\rm (g)}]   $B_4X=A_4$ has a solution $X \in {\mathbb C}^n_{\rm H}$ such that
$A_1-B_1XB_1^{*} \geqslant 0$ $(\, A_1-B_1XB_1^{*} \leqslant 0\,)$ if and only if
$r(M) = i_{+}(N)\ \ \left(r(M) = i_{-}(N) \right).$

\item[{\rm (f)}]  $A_1 - B_1XB_1^{*}  \geqslant 0$ $(\, A_1 - B_1XB_1^{*}  \geqslant 0\,)$  holds for all Hermitian solutions of $B_4X=A_4$ if and only if $i_{-}(N) = r(B_4)$ $(\, i_{+}(N) = r(B_4)\,).$
 \end{enumerate}
\end{theorem}

\noindent {\bf Proof.}  \ From Lemma \ref{T15b}(a), the general
Hermitian solution of $B_4X=A_4$ can be written as
\begin{align}
X=B_4^{\dag}A_4 + (B_4^{\dag}A_4)^{*} - B_4^{\dag}A_4B_4^{\dag}B_4 + F_{B_4}WF_{B_4},
\label{46}
\end{align}
where $W \in \mathbb C^n_{{\rm H}}$ is arbitrary. Substituting (\ref{45}) into $A_1-B_1XB_1^{*}$
gives
\begin{align}
A_1-B_1XB_1^{*} = P - B_1F_{B_4}WF_{B_4}B_1^{*},
 \label{47}
\end{align}
where $G =  A_1-B_1B_4^{\dag}A_4B_1^{*} - B_1(B_4^{\dag}A_4)^{*}B_1^{*} + B_1B_4^{\dag}A_4B_4^{\dag}B_4B_1^{*}$.
Applying (\ref{qq133})--(\ref{qq136}) to (\ref{46}) yields
 \begin{align}
\max_{X\in {\mathcal S}}r(\,A_1-B_1XB_1^{*}\,) & = \max_{W \in \mathbb C^n_{{\rm H}}}
 \!\!r(\,G-B_1F_{B_4}WF_{B_4}B_1^{*}\,)= r[\, G, \, B_1F_{B_4}\,],
 \label{48}
\\
\min_{X\in {\mathcal S}} r(\,A_1-B_1XB_1^{*}\,) & = \min_{W \in \mathbb C^n_{{\rm H}}}
\!\!r(\,G-B_1F_{B_4}WF_{B_4}B_1^{*}\,)= 2r[\, G,\,B_1F_{B_4}\,] - r\!\left[\!\!\begin{array}{cc}
 G & B_1F_{B_4}
\\
F_{B_4}B_1^{*} & 0 \end{array} \!\!\right]\!,
 \label{49}
\\
\max_{X\in {\mathcal S}}i_{\pm}(\,A_1-B_1XB_1^{*}\,) &
= \max_{W \in \mathbb C^n_{{\rm H}}} \!\!r(\,G-B_1F_{B_4}WF_{B_4}B_1^{*}\,)=
i_{\pm}\!\left[\!\!\begin{array}{cc}  G & B_1F_{B_4}
\\
F_{B_4}B_1^{*} & 0 \end{array} \!\!\right]\!,
 \label{410}
\\
\min_{X\in {\mathcal S}}
i_{\pm}(\,A_1-B_1XB_1^{*}\,) & = \min_{W \in \mathbb C^n_{{\rm H}}}
\!\!r(\,G - B_1F_{B_4}WF_{B_4}B_1^{*}\,)= r[\, G,\,B_1F_{B_4}\,]
- i_{\mp}\!\left[\!\! \begin{array}{cc}  G & B_1F_{B_4}
\\
F_{B_4}B_1^{*} & 0 \end{array} \!\!\right]\!.
 \label{411}
\end{align}
It is easy to verify that under $B_4B_4^{\dag}A_4= A_4$, the
equality $B_4(B_4^{\dag}A_4)^{*}= B_4A_4^{*}(B_4^{\dag})^{*} = A_4
B_4^{*}(B_4^{\dag})^{*} = A_4B_4^{\dag}B_4$ holds. In this case,
applying It is easy to verify by (1.17) and (1.25)  to (\ref{48})--(\ref{411})
and simplifying by elementary matrix operations and
 congruence matrix operations, we obtain
\begin{align}
r[\, G, \, B_1F_{B_4}\,] &  = r\!\left[\!\! \begin{array}{cc}
A_1 - B_1B_4^{\dag}A_4B_1^{*}-B_1(B_4^{\dag}A_4)^{*}B_1^{*} +
B_1B_4^{\dag}A_4B_4^{\dag}B_4B_1^{*} & B_1  \\ 0 & B_4 \end{array}
\!\!\right]-r(B_4) \nb
\\
& =r\!\left[\!\! \begin{array}{cc}  A_1 & B_1  \\
A_4B_1^{*}  + B_4(B_4^{\dag}A_4)^{*}B_1^{*} -A_4B_4^{\dag}B_4B_1^{*} & B_4
\end{array} \!\!\right] - r(B_4) \nb
\\
& = r\!\left[\!\! \begin{array}{cc}  A_1 & B_1  \\  A_4B_1^{*}& B_4
\end{array} \!\!\right]-r(B_4) =  r(M) - r(B_4),
\label{412}
\end{align}
\begin{align}
& i_{\pm}\!\left[\! \begin{array}{cc}  G & B_1F_{B_4} \\
F_{B_4}B_1^{*} & 0 \end{array} \!\right] \nb
\\
 & =i_{\pm}\!\left[\!\! \begin{array}{ccc}  A_1-B_1B_4^{\dag}A_4B_1^{*}-B_1(B_4^{\dag}A_4)^{*}B_1^{*} +
B_1B_4^{\dag}A_4B_4^{\dag}B_4B_1^{*} & B_1 & 0 \\
B_1^{*} & 0  & B_4^{*} \\ 0 & B_4 & 0 \end{array} \!\!\right] - r(B_4)
\nb
\\
& = i_{\pm}\!\left[\!\! \begin{array}{ccc}  A_1 & B_1 &
\frac{1}{2}B_1B_4^{\dag}A_4B_4^{*}
 + \frac{1}{2}B_1A_4^{*} - \frac{1}{2}B_1B_4^{\dag}A_4B_4^{*} \\
B_1^{*} & 0  & B_4^{*} \\
\frac{1}{2}A_4B_1^{*} + \frac{1}{2}B_4B_1(B_4^{\dag}A_4)^{*}B_1^{*} -
\frac{1}{2}A_4B_4^{\dag}B_4B_1^{*} & B_4 & 0 \end{array} \!\!\right]
\nb
\\
&  \ \ \ \  - r(B_4) \nb
\\
& = i_{\pm}\!\left[\!\! \begin{array}{ccc}  A_1 & B_1 & \frac{1}{2}B_1A_4^{*} \\
B_1^{*} & 0  & B_4^{*} \\
\frac{1}{2}A_4B_1^{*} & B_4 & 0 \end{array} \!\!\right] - r(B_4)
=i_{\pm}\!\left[\!\! \begin{array}{ccc} A_1 & B_1 & 0 \\ B_1^{*} & 0 &
B_4^{*} \\ 0 & B_4 &  -A_4B_4^{*} \end{array} \!\!\right] - r(B_4) =  i_{\pm}(N) - r(B_4).
\label{413}
\end{align}
Substituting (\ref{412}) and (\ref{413}) into
(\ref{48})--(\ref{411}) yields (\ref{42})--(\ref{45}).
 Applying Lemma \ref{T12} to  (\ref{42})--(\ref{45}) yields
(a)--(f). \qquad $\Box$

\begin{theorem} \label{T42}
Assume that the matrix equation $B_4X = A_4$ in {\rm (\ref{yy12})}
has a nonnegative definite solution$,$  i.e.$,$ ${\mathscr R}(A_4) \subseteq {\mathscr
R}(B_4),$ $A_4B_4^{*} \geqslant 0$ and $r(A_4B_4^{*}) =r(A_4),$
 and let
\begin{align}
{\mathcal S} = \{\, 0 \leqslant X\in {\mathbb C}^{n}_{{\rm H}} \ | \  A_4X = B_4 \, \},  \
M_1 = \left[\!\! \begin{array}{cc}  A_1 & B_1 \\
A_4B_1^{*} & B_4 \end{array} \!\!\right]\!,  \ M_2 =\left[\!\! \begin{array}{cc} A_1 & B_1A_4^{*} \\  A_4B_1^{*}  & A_4B_4^{*}  \end{array} \!\right], \  N =\left[\! \begin{array}{ccc}  A_1 & B_1 & 0
\\
B_1^{*} & 0 & B_4^{*} \\ 0 & B_4 &  -A_4B_4^{*}  \end{array} \!\!\right]\!.
\label{414}
\end{align}
Then$,$
 \begin{align}
\max_{X\in {\mathcal S}} r(\,A_1 - B_1XB_1^{*}\,) & = r(M_1) - r(B_4),
 \label{415}
\\
\min_{X\in {\mathcal S}} r(\,A_1-B_1XB_1^{*}\,) & = r(M_1) + i_{-}(M_2) - i_{-}(N),
 \label{416}
\\
\max_{X\in {\mathcal S}} i_{+}(\,A_1-B_1XB_1^{*}\,) & = i_{+}(M_2) - r(A_4),
 \label{417}
\\
\min_{X\in {\mathcal S}}i_{+}(\,A_1-B_1XB_1^{*}\,) & = r(M_1) - i_{-}(N),
 \label{418}
 \\
 \max_{X\in {\mathcal S}} i_{-}(\,A_1-B_1XB_1^{*}\,) & = i_{-}(N) - r(A_4),
 \label{419}
\\
\min_{X\in {\mathcal S}}i_{-}(\,A_1-B_1XB_1^{*}\,) & = i_{-}(M_2).
 \label{420}
\end{align}
In  consequences$,$
\begin{enumerate}
\item[{\rm (a)}] $B_4X=A_4$ has a nonnegative definite solution such that $A_1-B_1XB_1^{*}$  is nonsingular
 if and only if $r(M_1)=r(B_4)+m_1.$

\item[{\rm (b)}]  $A_1-B_1XB_1^{*}$  is nonsingular for all nonnegative definite solution of $B_4X=A_4$
 if and only if $r(M_1) + i_{-}(M_2) = i_{-}(N) + m_1.$

\item[{\rm (c)}] The pair of matrix equations  $B_1XB_1^{*} = A_1$ and $B_4X=A_4$ have a common nonnegative definite solution
 if and only if $r(M_1) + i_{-}(M_2) = i_{-}(N).$

\item[{\rm (d)}]  $B_1XB_1^{*} = A_1$ holds for all nonnegative definite solutions of $B_4X=A_4$
 if and only if $r(M) = r(B_4).$

\item[{\rm (e)}] $B_4X=A_4$ has a solution $X \in {\mathbb C}^n_{\rm H}$ such that $ A_1 - B_1XB_1^{*} >0$ if and only if $i_{+}(M_2) = r(A_4) + m_1.$

\item[{\rm (f)}]  $A_1 - B_1XB_1^{*} >0$  holds for all Hermitian solutions of $B_4X=A_4$
 if and only if $r(M_1)= i_{-}(N) + m_1.$

\item[{\rm (g)}] $B_4X=A_4$ has a solution $X \in {\mathbb C}^n_{\rm H}$ such that $ A_1 - B_1XB_1^{*} <0$ if and only if $i_{-}(N) = r(A_4) + m_1.$

\item[{\rm (h)}]  $A_1 - B_1XB_1^{*} <0$  holds for all Hermitian solutions of $B_4X=A_4$
 if and only if $i_{-}(M_2)= m_1.$

\item[{\rm (i)}] $B_4X=A_4$ has a solution $X \in {\mathbb C}^n_{\rm H}$ such that
$A_1-B_1XB_1^{*} \geqslant 0$  if and only if $M_2 \geqslant 0.$

\item[{\rm (j)}] $A_1 - B_1XB_1^{*} \geqslant 0$ holds for all Hermitian solutions of $B_4X=A_4$
 if and only if $i_{-}(N) = r(A_4).$

\item[{\rm (k)}] $B_4X=A_4$ has a solution $X \in {\mathbb C}^n_{\rm H}$ such that
$A_1-B_1XB_1^{*} \leqslant 0$  if and only if $r(M_1) = i_{-}(N).$

\item[{\rm (l)}] $A_1 - B_1XB_1^{*}  \leqslant 0$  holds for all Hermitian solutions of $B_4X=A_4$ if
and only if $i_{+}(M_2) = r(A_4).$
 \end{enumerate}
\end{theorem}

\noindent {\bf Proof.}  \ From Lemma \ref{T15b}(b), the general nonnegative definite solution of $B_4X=A_4$ can be written as
\begin{align}
X=A_4^{*}(A_4B_4^{*})^{\dag}A_4 + F_{B_4}WF_{B_4},
\label{421}
\end{align}
where $0 \leqslant W\in {\mathbb C}^{n}_{{\rm H}}$ is arbitrary. Substituting
(\ref{421}) into $A_1-B_1XB_1^{*}$ gives
\begin{align}
A_1-B_1XB_1^{*} = G - B_1F_{B_4}WF_{B_4}B_1^{*},
 \label{422}
\end{align}
where $G =  A_1- B_1 A_4^{*}(A_4B_4^{*})^{\dag}A_4B_1^{*}$.  Applying (\ref{qq140})--(\ref{qq142}) to (\ref{422}) yields
 \begin{align}
\max_{X\in {\mathcal S}} r(\,A_1-B_1XB_1^{*}\,) & = \max_{0 \leqslant W \in \mathbb C^n_{{\rm H}}}
 \!\! r(\,G - B_1F_{B_4}WF_{B_4}B_1^{*}\,)= r[\, G, \, B_1F_{B_4}\,],
 \label{423}
\\
\min_{X\in {\mathcal S}} r(\,A_1-B_1XB_1^{*}\,) & = \min_{0 \leqslant W \in \mathbb C^n_{{\rm H}}}
\!\!r(\,G - B_1F_{B_4}WF_{B_4}B_1^{*}\,) =  i_{-}(G) + r[\, G,\,B_1F_{B_4}\,] - i_{-}\!\left[\!\!\begin{array}{cc}  G & B_1F_{B_4} \\
F_{B_4}B_1^{*} & 0 \end{array} \!\!\right]\!,
 \label{424}
\\
\max_{X\in {\mathcal S}}i_{+}(\,A_1-B_1XB_1^{*}\,) & = \max_{0 \leqslant W \in \mathbb C^n_{{\rm H}}} \!\!r(\,G - B_1F_{B_4}WF_{B_4}B_1^{*}\,)
= i_{+}(G),
 \label{425}
\\
 \min_{X\in {\mathcal S}}i_{+}(\,A_1-B_1XB_1^{*}\,) & = \max_{0 \leqslant W \in \mathbb C^n_{{\rm H}}} \!\!i_{+}(\,G -B_1F_{B_4}WF_{B_4}B_1^{*}\,)
  = r[\, G,\,B_1F_{B_4}\,] - i_{-}\!\left[\!\!\begin{array}{cc}  G & B_1F_{B_4} \\ F_{B_4}B_1^{*} & 0
   \end{array} \!\!\right]\!,
\label{426}
\\
\max_{X\in {\mathcal S}}i_{-}(\,A_1-B_1XB_1^{*}\,) & = \max_{0 \leqslant W \in \mathbb C^n_{{\rm H}}} \!\!i_{-}(\,G - B_1F_{B_4}WF_{B_4}B_1^{*}\,)
= i_{-}\!\left[\!\!\begin{array}{cc}  G & B_1F_{B_4} \\
F_{B_4}B_1^{*} & 0 \end{array} \!\!\right]\!,
 \label{427}
\\
\min_{X\in {\mathcal S}} i_{\pm}(\,A_1-B_1XB_1^{*}\,) & =
\min_{0 \leqslant W \in \mathbb C^n_{{\rm H}}} \!\!i_{-}(\,G - B_1F_{B_4}WF_{B_4}B_1^{*}\,)= i_{-}(G).
 \label{428}
\end{align}
It is easy to verify by (1.17), (1.23)  and (1.25) that
\begin{align}
r[\, G, \, B_1F_{B_4}\,] &  = r\!\left[\!\! \begin{array}{cc}
 A_1- B_1 A_4^{*}(A_4B_4^{*})^{\dag}A_4B_1^{*} & B_1  \\ 0 & B_4 \end{array}
\!\!\right]-r(B_4)  = r\!\left[\!\! \begin{array}{cc}  A_1 & B_1  \\  A_4B_1^{*}& B_4
\end{array} \!\!\right]-r(B_4),
\label{429}
\\
i_{\pm}\!\left[\! \begin{array}{cc}  G & B_1F_{B_4} \\
F_{B_4}B_1^{*} & 0 \end{array} \!\right]  & =i_{\pm}\!\left[\!\! \begin{array}{ccc}  A_1- B_1 A_4^{*}(A_4B_4^{*})^{\dag}A_4B_1^{*} & B_1 & 0 \\
B_1^{*} & 0  & B_4^{*} \\ 0 & B_4 & 0 \end{array} \!\!\right] - r(B_4) \nb
\\
& = i_{\pm}\!\left[\!\! \begin{array}{ccc}  A_1 & B_1 & \frac{1}{2}B_1A_4^{*} \\
B_1^{*} & 0  & B_4^{*} \\
\frac{1}{2}A_4B_1^{*} & B_4 & 0 \end{array} \!\!\right] - r(B_4)
=i_{\pm}\!\left[\!\! \begin{array}{ccc} A_1 & B_1 & 0 \\ B_1^{*} & 0 &
B_4^{*} \\ 0 & B_4 &  -A_4B_4^{*} \end{array} \!\!\right] - r(B_4).
\label{430}
\\
i_{\pm}(G) & =  i_{\pm}[\, A_1- B_1 A_4^{*}(A_4B_4^{*})^{\dag}A_4B_1^{*} \,]
 = i_{\pm}\!\left[\! \begin{array}{cc} A_1 & B_1A_4^{*} \\  A_4B_1^{*}  & A_4B_4^{*}  \end{array} \!\right]- i_{\pm}(A_4B_4^{*}).
 \label{431}
\end{align}
Substituting (\ref{429})--(\ref{431}) into (\ref{423})--(\ref{428}) yields (\ref{415})--(\ref{420}).
Applying Lemma \ref{T12} to  (\ref{42})--(\ref{45}) yields
(a)--(l).  \qquad $\Box$

\begin{corollary} \label{T44}
Assume that the matrix equation in Lemma 1.10 has a Hermitian solution$,$ $P\in {\mathbb C}^{n}_{{\rm H}},$
and let ${\mathcal S} = \{\, X\in {\mathbb C}^{n}_{{\rm H}} \ | \  AX = B \, \}.$
Then$,$
\begin{align}
\max_{X \in {\mathcal S}}r(\, X - P\,) & = r(\,B -  AP \,) - r(A) + n,
\label{432}
\\
\min_{X \in {\mathcal S}}r(\, X - P \,) & =  2r(\, B - AP \,) - r(\,BA^{*} - APA^{*} \,),
\label{433}
\\
\max_{X \in {\mathcal S}} i_{\pm}(\, X - P \,) & = i_{\pm}(\,BA^{*} - APA^{*} \,) - r(A) + n,
\label{434}
\\
\min_{X \in {\mathcal S}} i_{\pm}(\, X - P \,) & =  r(\,B -  AP \,) - i_{\mp}(\,BA^{*} - APA^{*} \,).
\label{435}
\end{align}
In consequence$,$
\begin{enumerate}
\item[{\rm(a)}] There exists an $X \in {\mathcal S}$ such that $X - P$ is nonsingular
   if and only if
$$
{\mathscr R}(\, AP - B\,) = {\mathscr R}(A).
$$

\item[{\rm(b)}]  $X - P$ is nonsingular for all $X \in {\mathcal S}$ if and only if
$$
2r(\,B -  AP\,) = r(\,BA^{*} - APA^{*} \,) + n.
$$

\item[{\rm(c)}] There exists an $X \in {\mathcal S}$ such that $X > P$ $(X < P)$ holds
   if and only if
$$
{\mathscr R}(\, BA^{*} - APA^{*} \,) = {\mathscr R}(A) \ \ and \ \ BA^{*} \geqslant APA^{*}\ \
\left( \, {\mathscr R}(\, BA^{*} - APA^{*} \,) = {\mathscr R}(A) \ \ and \ \ BA^{*} \leqslant APA^{*} \,\right).
$$

\item[{\rm(d)}]  $X > P$ $(X < P)$ holds for all $X \in {\mathcal S}$ if and only if
$$
r(\,B -  AP\,) = n  \ \ and \ \  BA^{*} \geqslant APA^{*}\ \  \left(\,r(\, B - AP\,) = n  \ \ and \ \
 AB^{*} \leqslant APA^{*} \, \right).
$$

\item[{\rm(e)}] There exists an $X \in {\mathcal S}$ such that $X \geqslant P$
 $(X \leqslant P)$  holds
   if and only if
$$
 {\mathscr R}(\,B -  AP \,) = {\mathscr R}(\, BA^{*} - APA^{*} \,)  \ \ and  \ \
  BA^{*} \geqslant APA^{*} \ \  \left( \, {\mathscr R}(\, B - AP \,) = {\mathscr R}(\,BA^{*} - APA^{*}\,)  \ \ and  \ \
  BA^{*} \leqslant APA^{*}  \, \right)\!.
$$

\item[{\rm(f)}]  $X \geqslant P$ $(X \leqslant P)$  holds for all $X \in {\mathcal S}$ if and only if
$$
 BA^{*} \geqslant APA^{*} \ \ and  \ \   r(A) = n  \ \ \left(\, BA^{*} \leqslant APA^{*} \ \ and  \ \   r(A) = n  \, \right).
$$
\end{enumerate}
\end{corollary}

\begin{corollary} \label{T45}
Assume that the matrix equation in Lemma 1.10  has a Hermitian solution $X  \geqslant 0,$
and let $ 0 \leqslant  P\in {\mathbb C}^{n}_{{\rm H}}.$ Also$,$ define
\begin{align}
{\mathcal S} = \{\, 0 \leqslant X\in {\mathbb C}^{n}_{{\rm H}} \ | \  AX = B \, \}, \ \ M = \left[\! \begin{array}{cc}
 BA^{*} & B \\ B^{*} & P \end{array} \!\right]\!.
\label{436}
\end{align}
Then$,$
\begin{align}
\max_{X \in {\mathcal S}}r(\, X - P) & = r(\,B -  AP \,) - r(A) + n,
\label{437}
\\
\min_{X \in {\mathcal S}}r(\, X - P \,) & = i_{-}(M) + r(\,B -  AP \,) - i_{+}(\,BA^{*} - APA^{*} \,),
\label{438}
\\
\max_{X \in {\mathcal S}} i_{+}(\, X - P \,) &  = i_{+}(\,BA^{*} - APA^{*} \,) - r(A) + n,
\label{439}
\\
\min_{X \in {\mathcal S}} i_{+}(\, X - P \,) & = i_{-}(M),
\label{440}
\\
\max_{X \in {\mathcal S}} i_{-}(\, X - P \,) &  = i_{+}(M) - r(B),
\label{441}
\\
\min_{X \in {\mathcal S}} i_{-}(\, X - P \,) & =  r(\,B -  AP \,) - i_{+}(\,BA^{*} - APA^{*} \,).
\label{442}
\end{align}
In consequence$,$
\begin{enumerate}
\item[{\rm(a)}] There exists an $X \in {\mathcal S}$ such that $X - P$ is nonsingular
   if and only if ${\mathscr R}(\, B - AP \,) = {\mathscr R}(A).$

\item[{\rm(b)}]  $X - P$ is nonsingular for all $X \in {\mathcal S}$ if and only if
$i_{-}(M) +  r(\,B -  AP \,) = i_{+}(\,BA^{*} - APA^{*} \,) +n.$

\item[{\rm(c)}] There exists an $X \in {\mathcal S}$ such that $X > P$  holds
   if and only if ${\mathscr R}(\, BA^{*} - APA^{*} \,) = {\mathscr R}(A)$ and $BA^{*} \geqslant APA^{*}.$

\item[{\rm(d)}]  $X > P$ holds for all $X \in {\mathcal S}$ if and only if
$i_{-}(M) = r(A).$

\item[{\rm(e)}] There exists an $X \in {\mathcal S}$ such that $X < P$  holds
   if and only if $i_{-}(M) = r(B) + n.$

\item[{\rm(f)}]  $X < P$ holds for all $X \in {\mathcal S}$ if and only if
$r(\, B - AP \,) = n$ and $BA^{*} \leqslant APA^{*}.$

\item[{\rm(g)}] There exists an $X \in {\mathcal S}$ such that $X \geqslant P$  if and only if
${\mathscr R}(\,B -  AP \,) = {\mathscr R}(\, BA^{*} - APA^{*} \,)$  and  $BA^{*} \geqslant APA^{*}.$

\item[{\rm(h)}]  $X \geqslant P$  holds for all $X \in {\mathcal S}$ if and only if
$i_{-}(M) = r(B).$

\item[{\rm(i)}] There exists an $X \in {\mathcal S}$ such that $X \leqslant P$  if and only if
$M \geqslant 0.$

\item[{\rm(j)}] $X \leqslant P$  holds for all $X \in {\mathcal S}$ if and only if
$i_{+}(\,BA^{*} - APA^{*} \,) = n  - r(A).$
\end{enumerate}
\end{corollary}

\begin{corollary} \label{T46}
Assume that the matrix equation in Lemma 1.10  has a Hermitian solution$.$ Then$,$
 \begin{align}
\max_{AX=B,  \, X\in \mathbb C^n_{{\rm H}}} r(X) & = n + r(B) -
r(A),
 \label{443}
\\
\min_{AX=B, \, X\in \mathbb C^n_{{\rm H}}} r(X) & = 2r(B)- r(AB^{*}),
 \label{444}
 \\
 \max_{AX=B,  \, X\in \mathbb C^n_{{\rm H}}} i_{\pm}(X) & = n + i_{\pm}(AB^{*}) - r(A),
 \label{445}
\\
\min_{AX=B, \, X\in \mathbb C^n_{{\rm H}}} i_{\pm}(X) & = r(B)-
i_{\mp}(AB^{*}).
 \label{446}
\end{align}
Hence$,$
\begin{enumerate}
\item[{\rm (a)}] $AX=B$ has a nonsingular Hermitian solution if and only if  $r(A) = r(B).$

\item[{\rm (b)}]  $AX=B$ has a solution $X >0$ $(X <0)$ if and only if $AB^{*} \geqslant 0$  and $r(AB^{*}) = r(A)$
$(AB^{*} \leqslant 0$  and $r(AB^{*}) = r(A)).$

\item[{\rm (c)}]  $AX=B$ has a solution $X \geqslant 0$ $(X \leqslant 0)$ if and only if
$AB^{*} \geqslant 0$  and $r(AB^{*}) = r(B)$ $(AB^{*} \leqslant 0$  and
$r(AB^{*}) = r(B)).$

\item[{\rm (d)}] The rank of the Hermitian solution of $AX=B$ is invariant
 $\Leftrightarrow$ the positive index of  inertia of the Hermitian solution of $AX=B$ is invariant
  $\Leftrightarrow$ the negative index of inertia of the Hermitian solution of $AX=B$ is invariant
 $\Leftrightarrow$ $r(AB^{*}) = r(A) + r(B) - n.$
  \end{enumerate}
 \end{corollary}

Finally, we rewrite the matrix equation $AX=B$ as
\begin{align}
[\,A_1, \ A_2\,] \left[\!\!\begin{array}{cc}  X_1 & X_2 \\
X_2^{*} & X_3 \end{array} \!\!\right]=[\,B_1, \, B_2\,],
\label{447}
\end{align}
where $A_i \in {\mathbb C}^{m \times n_i}$,  $B_i \in {\mathbb C}^{m
\times n_i}$, $X_1 \in \mathbb C^{n_1}_{{\rm H}}$, $X_2 \in \mathbb
C^{n_1 \times n_2}$, $X_3 \in \mathbb C^{n_2}_{{\rm H}}$ for  $i=1,
\,2$ and $n_1+n_2=n.$ Note that the unknown submatrices in
(\ref{447}) can be  written as
\begin{align}
X_1=P_1XP_1^{*}, \ \ \  X_2=P_1XP_2^{*},  \ \ \ X_3=P_2XP_2^{*},
\label{448}
\end{align}
where $P_1=[\, I_{n_1}, \, 0\,]$ and $P_2=[\,0, \, I_{n_2}\,]$.  We
next find the extremal ranks and inertias of the submatrices $X_1$
and $X_3$ in a Hermitian solution of (\ref{447}). For convenience,
let
\begin{align}
{\cal S}_1=\{X_1 \in \mathbb C^{n_1}_{\rm H} \ | \ X_1=P_1XP_1^{*},\,
AX = B, \ X \in \mathbb C^{n}_{\rm H} \},
 \label{449}
 \\
{\cal S}_3=\{X_3 \in \mathbb C^{n_2}_{\rm H}  \ | \ X_3=P_2XP_2^{*},\,
AX=B,  \  X \in \mathbb C^{n}_{\rm H} \}.
\label{450}
\end{align}
Applying Theorem \ref{T41} to (\ref{449}) and (\ref{450}) gives the
following results. The details  of the proof are omitted.

\begin{theorem} \label{T47}
Assume that matrix equation in {\rm (\ref{417})} has a Hermitian
solution$,$ and let ${\cal S}_1$  and ${\cal S}_3$ be of the forms
in  {\rm (\ref{449})} and {\rm (\ref{450})}$.$ Then$,$ the global
maximal and minimal ranks and inertias of the Hermitian matrices in
${\cal S}_1$  and ${\cal S}_3$ are given by
 \begin{align}
\max_{X_1 \in {\cal S}_1} r(\,X_1\,) & = n_1 + r[\,A_2, \, B_1\,] -
r(A),
 \label{451}
\\
\min_{X_1 \in {\cal S}_1} r(\,X_1\,) & = 2r[\, A_2, \, B_1 \,] -
r\!\left[\!\! \begin{array}{cc}
 AB^{*} & A_2
\\
 A_2^{*} & 0 \end{array} \!\!\right]\!,
 \label{452}
\\
\max_{X_1 \in {\cal S}_1} i_{\pm}(\,X_1\,) & =  n_1 +
i_{\pm}\!\left[\!\! \begin{array}{cc}
 AB^{*} & A_2
\\
 A_2^{*} & 0 \end{array} \!\!\right] - r(A),
 \label{453}
\\
\min_{X_1 \in {\cal S}_1} i_{\pm}(\,X_1\,) & = r[\,A_2, \, B_1\,] -
i_{\mp}\!\left[\!\! \begin{array}{cc}
 AB^{*} & A_2
\\
 A_2^{*} & 0 \end{array} \!\!\right]\!,
 \label{454}
 \end{align}
 and
 \begin{align}
 \max_{X_3 \in {\cal S}_3} r(\,X_3\,) & = n_2 + r[\,A_1, \, B_2\,] - r(A),
 \label{455}
\\
\min_{X_3 \in {\cal S}_3} r(\,X_3\,) & = 2r[\, A_1, \, B_2 \,] -
r\!\left[\!\! \begin{array}{cc}
 AB^{*} & A_1
\\
 A_1^{*} & 0 \end{array} \!\!\right]\!,
 \label{456}
\\
\max_{X_3 \in {\cal S}_3} i_{\pm}(\,X_3\,) & =  n_2 +
i_{\pm}\!\left[\!\! \begin{array}{cc}
 AB^{*} & A_1
\\
 A_1^{*} & 0 \end{array} \!\!\right] - r(A),
 \label{457}
\\
\min_{X_3 \in {\cal S}_3} i_{\pm}(\,X_3\,) & = r[\,A_1, \, B_2\,] -
i_{\mp}\!\left[\!\! \begin{array}{cc}
 AB^{*} & A_1
\\
 A_1^{*} & 0 \end{array} \!\!\right]\!.
 \label{458}
\end{align}
\end{theorem}

Applying  Lemma \ref{T12} to  (\ref{451})--(\ref{454}),  we easily
obtain the following algebraic properties of the submatrix $X_1$ in
(\ref{447}).

\begin{corollary} \label{T48}
Assume that matrix equation in {\rm (\ref{447})} has a Hermitian
solution$.$ Then$,$
\begin{enumerate}
\item[{\rm(a)}] {\rm (\ref{447})} has a Hermitian solution in which $X_1$ is nonsingular if and
only if $r[\,A_2, \, B_1\,] = r(A).$

\item[{\rm(b)}]   $X_1$ is nonsingular in all Hermitian solutions of  {\rm (\ref{447})}
if and only if  $ r\!\left[\!\! \begin{array}{cc}
 AB^{*} & A_2 \\ A_2^{*} & 0 \end{array} \!\!\right] = 2r[\, A_2, \, B_1 \,] -n_1.$

\item[{\rm(c)}] {\rm (\ref{447})} has a Hermitian solution in which $X_1 >0$ $(X_1 <0)$ if and
only if
$$
i_{+}\!\left[\!\! \begin{array}{cc} AB^{*} & A_2 \\  A_2^{*} & 0
\end{array} \!\!\right] = r(A) \ \ \left(i_{-}\!\left[\!\!
\begin{array}{cc} AB^{*} & A_2 \\  A_2^{*} & 0 \end{array} \!\!\right] =
r(A)\right)\!.
$$

\item[{\rm(d)}]   $X_1 >0$ $(X_1 <0)$ in all Hermitian solutions of  {\rm (\ref{447})}
if and only if
$$
i_{-}\!\left[\!\! \begin{array}{cc}
 AB^{*} & A_2 \\ A_2^{*} & 0 \end{array} \!\!\right] = r[\,A_2, \, B_1\,] - n_1 \ \
 \left( i_{+}\!\left[\!\! \begin{array}{cc} AB^{*} & A_2 \\ A_2^{*} & 0 \end{array} \!\!\right]
 = r[\,A_2, \, B_1\,] - n_1\right)\!.
$$

\item[{\rm(e)}] {\rm (\ref{447})} has a Hermitian solution in which $X_1 \geqslant 0$ $(X_1 \leqslant 0)$
if and only if
$$
i_{+}\!\left[\!\! \begin{array}{cc} AB^{*} & A_2 \\  A_2^{*} & 0
\end{array} \!\!\right] = r[\,A_2, \, B_1\,]  \ \
\left(i_{-}\!\left[\!\! \begin{array}{cc} AB^{*} & A_2 \\  A_2^{*} & 0
\end{array} \!\!\right] = r[\,A_2, \, B_1\,]\right)\!.
$$

\item[{\rm(f)}] $X_1 \geqslant 0$ $(X_1 \leqslant 0)$ in all Hermitian solutions of  {\rm (\ref{447})}
if and only if
$$
i_{-}\!\left[\!\! \begin{array}{cc}  AB^{*} & A_2
\\
 A_2^{*} & 0 \end{array} \!\!\right]=  r(A) - n_1 \ \   \ \ \left(i_{+}\!\left[\!\! \begin{array}{cc}
 AB^{*} & A_2
\\
 A_2^{*} & 0 \end{array} \!\!\right]=  r(A) - n_1 \right)\!.
$$

\item[{\rm(g)}] {\rm (\ref{447})}  has a Hermitian  solution in which $X_1 = 0$ if and only if
${\mathscr R}(B_1)\subseteq {\mathscr R}(A_2).$

\item[{\rm(h)}] $X_1 =0 $ in all Hermitian solutions of  {\rm (\ref{447})}
if and only if  $r[\,A_2, \, B_1\,] = r(A) - n_1.$

\item[{\rm (i)}] The rank of $X_1$ in the Hermitian solution of {\rm (\ref{447})} is invariant
 $\Leftrightarrow$  the positive index of inertia of $X_1$ in the Hermitian solution of
 {\rm (\ref{447})} is invariant  $\Leftrightarrow$  the negative index of inertia of $X_1$
  in the Hermitian solution of {\rm (\ref{447})} is invariant
 $\Leftrightarrow$  $r\!\left[\!\! \begin{array}{cc}
 AB^{*} & A_2 \\ A_2^{*} & 0 \end{array} \!\!\right] = r[\,A_2, \, B_1\,] + r(A) - n_1.$

\end{enumerate}
\end{corollary}

\section{Conclusions}

\renewcommand{\theequation}{\thesection.\arabic{equation}}
\setcounter{section}{5} \setcounter{equation}{0}

In this paper, we studied the problems of maximizing and minimizing the rank and inertia
of the constrained matrix expression in (\ref{yy11}) and (\ref{yy12}), and obtained many
symbolic formulas for calculating the extremal ranks and inertias of (\ref{yy11}) and (\ref{yy12})
by using pure algebraic operations of matrices and their generalized inverses. As direct applications,
we gave necessary and sufficient conditions for the existence of $X$ satisfying the triple matrix
equations in (\ref{yy11}) and (\ref{yy12}), as well as some matrix inequalities.
Although the problems of maximizing and minimizing ranks and inertias of matrices are
generally regarded as NP-hard, the results presented in this previous sections
 as well as the papers \cite{LZGZ,LT1,LT2,LT3,LT4,T-laa10,T-ela10,T-ela11,T-med,TL} show that many
 closed-form formulas for calculating global extremal ranks inertias of some simpler matrix expressions
 can be established symbolically by using some pure algebraic operations of matrices, while these explicit formulas
 can be used to solve many fundamental problems in matrix theory, as mentioned in the
beginning of this paper. All the results obtained in these papers are brand-new, but easy to
understand within the scope of elementary linear algebra. This series of fruitful researches show
that for many basic or classic problems like solvability of matrix equations and matrix inequalities,
we are still able to establish a variety of innovative results by some new methods.

Motivated by the fruitful results and the analytical methods used in this paper,
 we mention some research problems for further consideration:
\begin{enumerate}
\item[{\rm (a)}] A challenging task is to give the
closed-form for the general common solution of
$B_2 XB^{*}_2 = A_2$ and $B_3 XB^{*}_3 = A_3$ that satisfies $X >0$ ($<0$, $\geqslant 0$, $\leqslant 0$),
which is equivalent to solving the inequalities
\begin{align}
X_0 +  VF_B  +  F_B V^{*} + F_{B_2}U F_{B_3} + F_{B_3}U^{*}F_{B_2} > 0  \, ( <0,  \  \geqslant 0, \
 \leqslant 0).
\end{align}
Moreover, give the extremal rank and partial inertia of $A_1 - B_1XB^{*}_1$
subject to $B_2 XB^{*}_2 = A_2$ and $B_3 XB^{*}_3 = A_3$ and $\geqslant 0$.

\item[{\rm (b)}] Give the extremal ranks and inertias of the LHMF $A_1 - B_1XB^{*}_1$
subject to the common Hermitian solution of the $k-1$ consistent
linear matrix equations
$$
[\,B_2XB^{*}_2, \ldots, B_kXB_k^{*}\,] =[\,A_2, \ldots, A_k\,],
$$
 and to establish necessary and sufficient condition for the set of matrix equations
 $$
[\,B_1XB^{*}_1, \ldots,
B_kXB_k^{*}\,] = [\,A_1, \ldots, A_k\,]
$$
 to have a common Hermitian solution, as well as a common nonnegative definite solution.

\item[{\rm (c)}]  Give the extremal rank and partial inertia of $A_1 - B_1XB^{*}_1$
subject to a linear matrix inequality $B_2XB^{*}_2 \geqslant A_2$. In such a case, it is necessary to first
give analytical expression for the general Hermitian solution of $B_2XB^{*}_2 \geqslant A_2$.

\item[{\rm (d)}] Give the extremal ranks and inertias of $A_1 - B_1XB^{*}_1$
subject to  $B_2X = A_2$ and $X \geqslant 0$.

\end{enumerate}

Since linear algebra is a successful theory with essential applications in most scientific fields,
the methods and results in matrix theory are prototypes of many concepts and content in other
advanced branches of mathematics. In particular,  matrix equations and matrix inequalities in the
 L\"owner partial ordering, as well as generalized inverses of matrices were sufficiently extended to
 their counterparts for operators in a Hilbert space, or elements in a ring with involution, and their
  algebraic properties were extensively studied in the literature.
   In most cases, the conclusions on the complex matrices and their counterparts in general algebraic settings
    are analogous. Also, note that the results in this paper are derived from
 ordinary algebraic operations of the given matrices and their  generalized inverses. Hence,
 it is no doubt that most of the conclusions in this paper can trivially be extended to
  the corresponding  equations and inequalities for linear operators on a Hilbert space or elements
  in a ring with involution.

\end{document}